\documentclass[11pt,a4paper]{amsart} 
\usepackage[english]{babel}
\usepackage{amssymb,xspace,enumerate,graphics,ifthen, dsfont}
\usepackage{amstext}
\theoremstyle{plain}
\usepackage{amsbsy,amsfonts,latexsym,enumerate}

\newtheorem{theorem}{Theorem}[section]
\newtheorem{lemma}[theorem]{Lemma}
\newtheorem{fact}[theorem]{Fact}
\newtheorem{proposition}[theorem]{Proposition}
\newtheorem{corollary}[theorem]{Corollary}

\newtheorem{definition}[theorem]{Definition}
\newcommand{\pr}{\noindent {\it Proof~: }}

\def\leqs{\lesssim}
\def\geqs{\gtrsim}
\def\eqs{\eqsim}
\def\loc{\mathrm{loc}}

\def\cc{\mathbb{C}}
\def\rr{\mathbb{R}}
\def\nn{\mathbb{N}}

\newcommand{\diffp}[2]{\frac{\partial #1}{\partial #2}}
\newcommand{\mlabel}[1]{\label {#1}}
\newcommand{\cp}[2]{{\mathcal{P}}_{#1}(#2)}
\newcommand{\cq}[2]{{\mathcal{Q}}_{#1}(#2)}
\newcommand{\ch}[1]{{\mathcal{H}}^{#1}(D)}

\renewcommand{\over}[2]{\genfrac{}{}{0pt}{}{#1}{#2}}
\newcommand{\vol}{\mathrm{Vol}}

\renewcommand{\d}{\mathrm{d}}

\newcommand{\dz}{{|r(z)|}}
\newcommand{\dzeta}{|r(\zeta)|}

\newcommand{\hl}{{h_\Lambda(z,t)}}

\renewcommand{\ij}{(I)_j}
\newcommand{\ijkl}{(I)_{j,k,l}}
\newcommand{\ii}{(II)}
\newcommand{\iijk}{(II)_{j,k}}
\newcommand{\iiij}{(III)_j}
\newcommand{\iiijkl}{(III)_{j,k,l}}

\title{$\mathcal{H}^p$-corona problem and convex domains of finite type }
\author{William ALEXANDRE}
\address[W. Alexandre]{Univ. Lille, CNRS, UMR 8524 - Laboratoire Paul
Painlev\'e, F-59000 Lille, France}
\email{ william.alexandre@univ-lille.fr}

\begin{document}

\begin{abstract}
We prove that the $\mathcal{H}^p$-corona problem has a solution for convex domains of finite type in $\cc^n$, $n\geq 2$.
\end{abstract}
\maketitle
\section{Introduction and main result}
Let $D$ be a domain in $\cc^n$ and $f_1,\ldots, f_k$ be $k$ functions in $\ch\infty$, the algebra of bounded holomorphic functions on $D$. Assume that for all $z\in D$, the following inequality holds true for some $\delta>0$
\[
\sum_{j=1}^k |f_j(z)|^2\geq \delta^2.
\]
To solve the $\mathcal{H}^\infty$-Corona Problem on $D$ is then to find $k$ functions $g_1,\ldots, g_k$ in $\ch{\infty}$ such that
for all $z\in D$, 
\[
\sum_{j=1}^k g_j(z)f_j(z)=1.
\]
The $\mathcal{H}^\infty$-Corona Problem is solved by Carleson in \cite{Carl} when $D$ is the unit disc of $\cc$ but is still an open question when $n\geq 2$, even if $D$ is the ball or the polydisc. On the other side, Sibony in \cite{Sibo} and Forn{\ae}ss and Sibony in \cite{FornSibo}  construct bounded pseudoconvex domains with smooth boundary and data $f_1,\ldots, f_k$, such that the Corona Problem has no solution. 
It is an interesting question to know for which domains in $\cc^n$ the Corona Problem may have a solution. As pointed out by Amar in \cite{Amar}, being able to solve the $\mathcal{H}^p$-Corona Problem is a necessary condition to solve the $\mathcal{H}^\infty$-Corona Problem. Let us state the $\mathcal{H}^p$-Corona Problem.

\par\medskip 

We write \(D\) as the set \(D=\{ z\in\cc^n,\ r(z)<0\}\) where $r$ is a smooth function on $\cc^n$ such  that $\d r\neq 0$ on the boundary of $D$. For $\varepsilon\in\rr$, we denote by $bD_\varepsilon$ the boundary of $D_ \varepsilon:=\{ z\in\cc^n,\ r(z)<\varepsilon\}$,  and  by $\d\sigma_\varepsilon $ the euclidean area measure on $bD_\varepsilon $.

The Hardy space $\ch p$, $p>0$, is the set of holomorphic functions $f$ on $D$ such that
\[\|f\|_{\ch p}=\left(\sup_{\varepsilon >0}\int_{bD_{-\varepsilon} }|f(z)|^p\d\sigma_{-\varepsilon} (z)\right)^{\frac1p}<+\infty.\]
By passing to the (almost everywhere) radial limit function, we may see the  space $\ch p$ as a closed subspace of $L^p(bD)$ (see \cite{Kra}). 

To solve the $\mathcal{H}^p$-Corona Problem is to find for any $h\in\ch p$, $k$ functions $h_1,\ldots, h_k\in \ch p$ such that 

\begin{align*}
  \sum_{j=1}^kh_jf_j=h. 
\end{align*}

Amar solves in \cite{Amar} the $\mathcal{H}^p$-Corona Problem on the ball of $\cc^n$, $n\geq 2$, and for two generators (i.e. \(k=2\)), for all $1<p<\infty$. It is also solved by Andersson and Carlsson for 2 generators in \cite{AndeCarl1} and for $k$ generators in \cite{AndeCarl2,AndeCarl3} on strictly pseudoconvex domains. In \cite{Lin}, Lin proves that the $\mathcal{H}^p$-Corona Problem has a solution for $k$ generators, $k\geq2$, on the polydisc of $\cc^n$, $n\geq 2$, $1<p<+\infty$.

In this article, we solve the $\mathcal{H}^p$-Corona Problem for 2 generators on convex domains of finite type.
\begin{theorem}\mlabel{courrone}
 Let \(D\subset \cc^n\), \(n\geq 2\), be a bounded convex domains of finite type with smooth boundary. Let $f_1,f_2$ in $\ch \infty$ and $\delta$ in $\rr$ be such that $|f_1|^2+|f_2|^2\geq \delta^2>0$ on $D$.  Then for all \(1<p<\infty\), all \(h\in\ch p\), there exist \(h_1,\ h_2\in\ch p\) such that
 \(h=h_1f_1+h_2f_2\).
\end{theorem}
In order to establish Theorem \ref{courrone}, as  Amar \cite{Amar} and Andersson-Carlsson \cite{AndeCarl1} do, we follow Wolff's proof of the one variable Corona Theorem. We first put 
\begin{align*}
 g_i&=\frac{\overline{f_i}}{|f_1|^2+|f_2|^2},\quad i=1,2,\\
 \omega&=\frac{\overline{f_1}\,\overline\partial \, \overline{f_2} - \overline{f_2}\,\overline \partial\, \overline{f_1}}{(|f_1|^2+|f_2|^2)^2}.
\end{align*}
It follows that
$$f_1g_1+f_2g_2=1,\qquad \overline\partial g_1=-f_2\omega, \qquad \overline\partial g_2=f_1\omega \qquad \text{and} \qquad \overline\partial \omega=0.$$ 

For $h$ holomorphic in $D$, we have $\overline\partial (h\omega)=0$. So one can find $u$ such that $\overline \partial u =h\omega$. Setting
\begin{align*}
 h_1&=hg_1 +u f_2, \\
 h_2&=hg_2 -u f_1,
\end{align*}
we have
\begin{align*}
 h&=f_1h_1+f_2h_2,\\
 \overline\partial h_1
 &=h\overline \partial g_1 +f_2\omega h=0,\\
  \overline\partial h_2
 &=h\overline \partial g_2 -f_1\omega h=0.
 \end{align*}
Moreover, since $g_1$, $g_2$, $f_1$ and $f_2$ are bounded on $D$, if $h$ belongs to $\ch p$, $h_1$ and $h_2$ will also be in $\ch p$ provided that $u$ belongs to $L^p(bD)$.
So the proof of Theorem \ref{courrone} is reduced to find $u\in L^p(bD)$ such that $\overline\partial u=h\omega$, i.e. to solve a $\overline\partial$-equation with boundary estimates. 
\par\medskip
As in \cite{Carl}, \cite{Amar}, \cite{AndeCarl1}, Carleson measures are in the present paper  an essential tool in order to the solve the $\mathcal{H}^p$-Corona Problem.
They are defined using the homogeneous structure of the boundary of the domain. For convex domains, one should use McNeal polydiscs defined in \cite{McN, McN1, McN2}. Since we need many objects in order to define them, we postpone the definition of the set of Carleson measures $W^1(D)$ and the set of $(p,q)$-Carleson currents $W_{p,q}^1(D)$ to Section \ref{notation}. We also define, in Section \ref{notation}, $BMO(bD)$, the space of functions of bounded mean oscillation on $bD$. We denote by $W^0(D)$ the set of bounded measures and by $W^0_{p,q}(D)$ the set of $(p,q)$-currents with bounded measure coefficients. Then, for $\alpha\in]0,1[$, $W^{\alpha}(D)$ is the complex interpolate space $[W^{0}(D),W^1(D)]_\alpha$ and $W^\alpha_{p,q}=[W^{0}_{p,q}(D),W^1_{p,q}(D)]_\alpha$. 

If $h\omega$ were in $W^{1-\frac1p}_{0,1}(D)$, the existence of $u$ and thus Theorem \ref{courrone} would be a consequence of Theorem 2.10 of \cite{Ale2}. But in general, this is not the case. One has to construct more elaborate functions $g_1$, $g_2$ and $\omega$ such that $h\omega$ belongs to  $W^{1-\frac1p}_{0,1}(D)$. This construction is done by Carleson in the one variable case, but it seems too difficult to carry it out in several variables. Instead, we use Wolff's approch who notices that $|r||\omega|^2$ and $|r||\partial \omega|$ are Carleson measures in the $1$-dimensional case.

We use results of Jasiczak \cite{Jas} in order to prove (see Subsection \ref{subsection31}) that $\partial(h\omega)$ satisfies the hypothesis of Theorem \ref{ddbar} below. In this theorem and in the sequel, $A\leqs B$ means there exists a constant $c>0$ such that $A\leq cB$ and $A\eqs B$ that $A\leqs B$ and $B\leqs A$ both hold.
\begin{theorem}\label{ddbar}
 Let $D$ be a bounded convex domain of finite type in $\cc^n$, let $\theta$ be a $\d$-closed $(1,1)$-current. Then
 \begin{enumerate}[(i)]
  \item if $ |r|\, \theta$ belongs to $W^{1}_{1,1}(D)$, then there exists $v$ such that $\theta=i\partial\overline\partial v$ and $\|v\|_{BMO(bD)}\leqs \||r|\,\theta\|_{W_{1,1}^1(D)}$, uniformly with respect to $\theta$.
  \item if $|r|\, \theta$ belongs to $W^{1-\frac1p}_{1,1}(D)$, $1\leq p<+\infty$, then there exists $v$ such that $\theta=i\partial\overline\partial v$ and $\|v\|_{L^p(bD)}\leqs \||r|\,\theta\|_{W_{1,1}^{1-\frac1p} (D)}$,  uniformly with respect to $\theta$.
 \end{enumerate}
\end{theorem}
So Theorem \ref{ddbar} gives us a function $v$ such that $\partial\overline{\partial} v=\partial (h\omega)$. Since 
\begin{align*}
 \partial (\overline\partial v-h\omega) &= \partial\overline\partial v-\partial (h\omega)=0,\\
 \overline{\partial}(\overline\partial v-h\omega) &= \overline\partial^2v-\overline \partial (h\omega)=0,
\end{align*}
the $1$-form $\overline\partial v-h\omega$ is $\d$-closed and we can solve the $\d$-equation for $\overline\partial v-h\omega$ : there exists a unique function $w$ such that $\d \overline w=\overline\partial v-h\omega$ and $\overline w(0)=v(0)$, where  $0$ is any point in $D$. Since, $\overline\partial v-h\omega$ is a $(0,1)$-form, we have
$$\partial \overline w= 0\quad \text{ and } \quad \overline\partial \overline w= \overline\partial v-h\omega.$$
Therefore $w$ is holomorphic and $u=v-\overline w$ satisfies
\begin{align*}
 \overline\partial u
 &= \overline \partial v -\overline\partial\,\overline w\\
 &=h\omega.
\end{align*}
Moreover, since $v$ already belongs to $L^p(bD)$, $u$ is in $L^p(bD)$ if and only if $w$ belongs to  $\ch p$. We will prove that this is indeed the case in Subsection \ref{subsection32} by methods similar to Amar and Andersson-Carlsson's method. This will solve the $\mathcal{H}^p$-Corona Problem for 2 generators in  convex domains of finite type.
\par\medskip
Many of the proofs in the present paper rely on interpolation, in particular between $\ch 1$ and $BMOA(D)$. We thus have to know what are the intermediate spaces between them. We will prove in Section \ref{S-interpolation} that when $D$ is a convex domain of finite type, $[\ch 1,BMOA(D)]_{1-\frac1p}=\ch p$, $1<p<+\infty$.  This result is also true when $D$ is strictly pseudoconvex and we prove it in the same way. However, the proof requires some regularity conditions on the tools, in our case the $\varepsilon$-extremal basis, which define the homogeneous structure of the boundary of the domains,
itself used to define $BMO(bD)$ and $BMOA(D)$. 
When $D$ is strictly pseudoconvex, the basis used are smooth, but not McNeal's $\varepsilon$-extremal basis. We overcome this difficulty by using the Bergman metric. I want to thank Éric Amar and Pierre Portal for helping me to understand the proof of the interpolation between $\ch 1$ and $BMOA(D)$ when $D$ is strictly pseudoconvex.
\par\medskip
The article is thus organised as follows : in Section \ref{notation}, we introduce the tools and objects relative to the structure of homogeneous spaces of $D$. In Section \ref{S-interpolation}, we prove the interpolation results we need. In Section \ref{bazar}, we prove that $h\omega$ satisfies the hypothesis of Theorem \ref{ddbar} and that $w$ belongs to $\ch p$. In Section \ref{section5}, we prove Theorem \ref{ddbar}.

\section{Notations}\mlabel{notation}

For $z$ near $bD$, $\varepsilon>0$ and $v\in\cc^n$, $v\neq 0$,  we denote by 
$\tau(z,v,\varepsilon)$ the distance from $z$ to $\{r=r(z)+\varepsilon\}$ in the complex direction $v$ :  
\begin{align*}
\tau(z,v,\varepsilon)&:=\sup\{t>0, |r(z+\lambda v)-r(z)|<\varepsilon, \forall \lambda \in\cc,\ |\lambda |<t\},
\end{align*}
 Using these distances, we define $\varepsilon$-extremal basis $w_1^*,\ldots, w_n^*$ at the point $z$, as given in \cite{BCD} :  $w_1^* =\eta_z$ is the outer unit normal to $bD_{r(z)}$ at $z$ and if $w^*_1,\ldots, w^*_{i-1}$ are already defined, then $w_i^*$ is a unit vector orthogonal to $w_1^*,\ldots, w^*_{i-1}$ such that $\tau(z,w^*_i,\varepsilon)= \sup_{\over{v\perp w^*_1,\ldots, w^*_{i-1}}{\|v\|=1}} \tau(z,v,\varepsilon)$. When $D$ is strictly convex, $w_1^*$ is the outer unit normal to $bD_{r(z)}$ and we may choose any basis of $T^\cc_z bD_{r(z)}$ for  $w^*_2,\ldots, w^*_n$. Therefore, when $D$ is strictly convex, an $\varepsilon$-extremal basis at $z$ can be chosen smoothly depending on the point $z$. Unfortunately, this is not the case for convex domains of finite type (see \cite{Hef2}).\\
We put $\tau_i(z,\varepsilon)=\tau(z,w_i^*,\varepsilon)$, for $i=1,\ldots,n$.
We have for a strictly convex domain $\tau_1(z,\varepsilon)\eqs \varepsilon$ and $\tau_j(z,\varepsilon)\eqs\varepsilon^{\frac12}$ for $j=2,\ldots,n$. For a convex domain of finite type $m$, we only have $\varepsilon^{\frac12}\leqs \tau_n(z,\varepsilon)\leq\ldots\leq \tau_2(z,\varepsilon)\leqs \varepsilon^{\frac 1m}$,  uniformly with respect to $z$ and $\varepsilon$.
The McNeal polydisc centered at $z$ of radius $\varepsilon $ is the set $$\cp\varepsilon z:=\left\{ \zeta=z+\sum_{i=1}^n \zeta^*_iw_i^*\in\cc^n, |\zeta^*_i|<\tau_i(z,\varepsilon),\:i=1,\ldots, n\right\}.$$ 
McNeal's polydiscs are used in order to define a pseudodistance $\delta$. We set for $\zeta,z$ near $bD$
$$\delta(z,\zeta):=\inf\{\varepsilon >0,\ \zeta\in \cp \varepsilon z\}.$$

\begin{definition}
We say that a positive finite measure $\mu$ on $D$ is a {\em Carleson measure} and we write $\mu\in W^1(D)$ if 
$$\|\mu\|_{W^1(D)}:= \sup_{\over{z\in bD}{\varepsilon>0}}\frac {\mu(\cp \varepsilon z\cap D) }{\sigma(\cp \varepsilon z\cap bD )}<\infty.$$
\end{definition}
Now we defined the notion of Carleson current already used in \cite{Ale2} and \cite{Ale3}.
For $z\in\cc^n$ and $v$ a non zero vector we 
set (see \cite{BCD})
\begin{align*}
k(z,v)&:=\frac{|r(z)| }{\tau(z,v,|r(z)| )}.
\end{align*}
For a fixed $z$, the convexity of $D$ implies that the function defined by $v\mapsto k(z,v)$ if $v \neq 0$, $0$ otherwise, is a kind of non-isotropic norm which will play for us the role of weight in the definition of Carleson currents.
\begin{definition}
We say that a $(p,q)$-current $\mu$ of order $0$ with measure coefficients  is a {\em $(p,q)$-Carleson current} if
$$\|\mu\|_{W^1_{p,q}}:=\sup_{u_1,\ldots,u_{p+q}} \left\|\frac1{k(\cdot,u_1)\ldots k(\cdot,u_{p+q})}\left|\mu(\cdot)[u_1,\ldots,u_{p+q}]\right|\right\|_{W^1}<\infty,$$
where the supremum is taken over all smooth vector fields $u_1,\ldots, u_{p+q}$ which never vanish and where $|\mu(\cdot)[u_1,\ldots,u_{p+q}]|$ is the absolute value of the measure $\mu(\cdot)[u_1,\ldots,u_{p+q}]$.

We denote by $W^1_{p,q}(D)$ the set of all $(p,q)$-Carleson currents.
\end{definition}
Let $W^0$ be the set of positive bounded measures on $D$. For $\mu\in W^0$, we put $\|\mu\|_{W^0}:=\mu(D)$. Analogously to $W^1_{p,q}(D)$ we define $W^0_{p,q}(D)$:
\begin{definition}
We say that $\mu$ is a $(p,q)$-current with bounded measure coefficients and we write $\mu\in W^0_{p,q}(D)$ if
$$\|\mu\|_{W^0_{p,q}}:=\sup_{u_1,\ldots,u_{p+q}} \left\|\frac1{k(\cdot,u_1)\ldots k(\cdot,u_{p+q})}\left|\mu(\cdot)[u_1,\ldots,u_{p+q}]\right|\right\|_{W^0}<\infty,$$
where the supremum is taken over all smooth vector fields $u_1,\ldots, u_{p+q}$ which never vanish and where $|\mu(\cdot)[u_1,\ldots,u_{p+q}]|$ is the absolute value of the measure $\mu(\cdot)[u_1,\ldots,u_{p+q}]$.
\end{definition}
For all $\alpha\in ]0,1[$ the space $W^\alpha_{p,q}(D)$ will denote the complex interpolate space between $W^0_{p,q}(D)$ and $W^1_{p,q}(D)$. One can ``understand'' these spaces by the work  of  Amar and Bonami who proved in \cite{AB}
\begin{proposition}\label{propAB}
 A measure $\mu$ belongs to $W^\alpha(D)$, $\alpha\in ]0,1[$, if and only if there exists a Carleson measure $\mu_1$ and $f\in L^{\frac1{1-\alpha}}(bD, d\mu_1)$ such that $\mu=fd\mu_1$.

\end{proposition}

\section{The interpolation space $[\ch 1,BMOA(D)]_\theta$}\mlabel{S-interpolation}
Let us first define the spaces $BMO(bD)$ and $BMOA(D)$. For $f\in L^1_{\loc}(bD)$, we set 
\begin{align*}
 \|f\|_{BMO(bD)}= \sup_{z\in bD, r>0} \frac1{\vol(bD\cap\cp r z)} \int_{bD\cap\cp r z}\left|f(\zeta) - f_{bD\cap\cp r z}\right|\d\sigma(\zeta)
\end{align*}
where, for $U\subset bD$, $f_{U}=\frac1{\vol({U})}\int_Uf(\zeta) \d\sigma(\zeta)$ and $\vol({U})$ is the euclidean volume of $U$.
Reminding that $\ch 1$ is a closed subset of $L^1(bD)$, the space $BMOA(D)$ is the set :
$$BMOA(D)=\{f\in \ch 1\ / \ \|f\|_{BMO(bD)}<+\infty\}.$$
It is well known that $\|\cdot\|_{BMO(bD)}$ is not a norm because $\|f\|_{BMO(bD)}=0$ if $f$ is constant. Therefore, as in \cite{K-Li1}, we equip $BMOA(D)$ with the following norm defined  for $f\in BMOA(D)$ by
$$\|f\|_{BMOA(D)}=\|f\|_1+\|f\|_{BMO(D)}.$$
\par\medskip

Let us recall the definition of the interpolation space $[\ch p,BMOA(D)]_{\theta}$, $p\geq 1$, $\theta\in]0,1[$. First we equip $\ch p+BMOA(D)=\{\phi_0+\phi_1\ / \  \phi_0\in\ch p, \phi_1\in BMOA(D)\}$ with the norm : $\|\phi\|_{\ch p+ BMOA(D)}=\inf\{\|\phi_0\|_{\ch p} + \| \phi_1\|_{BMOA(D)} \ / \ \phi_0\in\ch p, \phi_1\in BMOA(D), \phi=\phi_0+\phi_1\}$. Then an element $f$ of $[\ch p,BMOA(D)]_{\theta}$ is a complex valued function such that there exists an application $\Phi:\{z\in\cc\ / \ 0\leq \Re z\leq 1\}\to \ch p+BMOA(D)$ which satisfies 
\begin{enumerate}[(i)]
 \item \label{condi} $f=\Phi(\theta)$,
 \item $\Phi$ is continuous,
 \item $\Phi$ is analytic on $\{z\in\cc\ / \ 0< \Re z< 1\}$ 
 \item $t\mapsto \Phi(it)$ and $t\mapsto \Phi(1+it)$ are continuous from $\rr$ to $\ch p$ and $BMOA(D)$ respectively, 
 \item \label{condv}$\lim_{|t|\to +\infty} \|\Phi(it)\|_{\ch p}=0$ and $\lim_{|t|\to +\infty} \|\Phi(1+it)\|_{BMOA(D)}=0$.
\end{enumerate}
The norm of $f\in [\ch p,BMOA(D)]_{\theta}$ is 
 $$\|f\|_{[\ch p,BMOA(D)]_{\theta}}=\inf_{\Phi} \max\left(\sup_{\rr}\|\Phi(it)\|_{\ch p},\sup_{\rr}\|\Phi(1+it)\|_{BMOA(D)}\right)$$
where the infimum  is taken over all $\Phi$ satisfying (\ref{condi}-\ref{condv}) (see \cite{BL}). We will prove in this section the following result :
\begin{theorem}\mlabel{interpolation}
Let $D$ be a convex domain of finite type and $q$ be in $]1,+\infty[$. Then $\ch q=[\ch 1,BMOA(D)]_{1-\frac1q}$ with equivalent norms.
\end{theorem}

We prove Theorem \ref{interpolation} by showing that $[\ch p,BMOA(D)]_{1-\frac pq}=\ch q$  when $1<p<q<+\infty$ and by extending this result to the case $p=1$ using Wolff's note \cite{Wol}.

\begin{lemma}\mlabel{first_injection}
Let $D$ be a convex domain of finite type, $1<p< q<+\infty$.
Then $\ch q\subset [\ch p ,BMOA(D)]_{1-\frac pq}$ and for $f\in \ch q$, $\|f\|_{[\ch p ,BMOA(D)]_{1-\frac pq}}\leqs \|f\|_{\ch q}$, uniformly with respect to $f$.
\end{lemma}

\pr Let $f$ be an element of $\ch q$ and let $\tilde f\in L^q(bD)$ be its boundary value. Since $L^q(bD)=[L^p(bD),L^\infty(bD)]_{1-\frac pq}$, there exists 
$\Psi:\{z\in\cc\ / \ 0\leq \Re z\leq 1\}\to L^p(bD)+L^\infty(bD)$ such that $\tilde f =\Psi({1-\frac pq})$, $\Psi$ is analytic on $\{z\in\cc\ / \ 0< \Re z< 1\}$, $t\mapsto \Psi(it)$ is continuous from $\rr$ to $L^p(bD)$, $t\mapsto \Psi(1+it)$ is continuous from $\rr$ to $L^\infty(bD)$ and both tends to $0$ when $|t|$ goes to $+\infty$.

Now let $S$ be the Szegö projector (see \cite{Kra}). The Szegö projector in linear thus $\Phi=S\circ\Psi$ is holomorphic on $\{z\in\cc\ / \ 0< \Re z< 1\}$. From \cite{Mc-S1}, Theorem 3.4 and 5.1, $S:L^p(bD)\to \ch p$ is continuous for all $1<p<+\infty$ and from \cite{K-Li1}, Theorem 5.6, $S:L^\infty(bD)\to BMOA(D)$ is also continuous. Therefore $\Phi$ is continuous,  
$t\mapsto \Phi(it)$ and $t\mapsto \Phi(1+it)$ are continuous from $\rr$ to $\ch p$ and $BMOA(D)$ respectively, $\lim_{|t|\to +\infty} \|\Phi(it)\|_{L^p(bD)}=0$  and $\lim_{|t|\to +\infty} \|\Phi(1+it)\|_{L^\infty(bD)}=0$. Morevover, since $f$ is already holomorphic and since $\tilde f$ is the boundary value of $f$, $f=S(\tilde f)=\Phi({1-\frac pq})$. Thus $f$ belongs to $[\ch p,BMOA(D)]_{1-\frac pq}$. 
Moreover, the continuity of $S$ gives, uniformly with respect to $\Psi$,
\begin{align*}
\|f\|_{[\ch p ,BMOA(D)]_{1-\frac pq}}
&\leq \max\left(\sup_{t\in\rr}\|S\circ \Psi(it)\|_{\ch p},\sup_{t\in\rr}\|S\circ \Psi(1+it)\|_{BMOA(D)}\right)\\
&\leqs \max\left(\sup_{t\in\rr}\|\Psi(it)\|_{L^p(bD)},\sup_{t\in\rr}\|\Psi(1+it)\|_{L^\infty(bD)}\right).
\end{align*}
Taking the infimum among all $\Psi$, we get
$$\|f\|_{[\ch p ,BMOA(D)]_{1-\frac pq}}\leqs\|\tilde f\|_{[L^p(bD) ,L^\infty(bD)]_{1-\frac pq}}= \|\tilde f\|_{L^q(bD)}=\|f\|_{\ch q}.$$
\qed
\par\medskip
We need to prove the converse inclusion which is more involved. We consider the following maximal functions. For $f\in L^1_{loc}(bD)$ and $z\in bD$, we set 
\begin{align*}
f^\#(z)&=\sup_{\zeta\in bD,\ \varepsilon>0\ / \ z\in\cp\varepsilon \zeta }\frac1{\vol(\cp\varepsilon \zeta \cap bD)}\int_{\cp\varepsilon \zeta \cap bD} |f-f_{\cp\varepsilon \zeta \cap bD}|\d\sigma.
\end{align*}
We aim at proving that $^\#:[\ch p,BMOA(D)]_{1-\frac pq} \to L^q(bD)$ is continuous when $1<p< q<+\infty$.
In order to establish the continuity of $^\#$, we introduce the maximal function $M$. For $f\in L^1_{loc}(bD)$ and $z\in bD$, we set 
\begin{align*}
M(f)(z)&=\sup_{\zeta\in bD,\ \varepsilon>0\ / \ z\in\cp\varepsilon \zeta }\frac1{\vol(\cp\varepsilon \zeta \cap bD)}\int_{\cp\varepsilon \zeta \cap bD} |f|\d\sigma. 
\end{align*}
We control $\|f\|_{L^q(bD)}$ by $\|f^\#\|_{L^q(bD)}$ with the following lemma :
\begin{lemma}\mlabel{Lemma2}
Let $1<q<+\infty$. The following inequality holds uniformly  with respect to $f\in L^q(bD)$ : $\|f\|_{L^q(bD)}\leq \|f^\#\|_{L^q(bD)}+\frac 1{\vol(bD)}\left|\int_{bD} f\d\sigma\right|$
\end{lemma}
\pr 
Lebesgue's differentiation theorem implies that $|f|\leq Mf $ almost everywhere so 
\begin{align}
\|f\|_{L^p(bD)}\leq \|Mf\|_{L^p(bD)}. \mlabel{eq4} 
\end{align}
From \cite{Bur}, Theorem 2, for all $g\in L^1_{loc}(bD)$ such that $\int_{bD} g\d\sigma=0$, we have $\|Mg \|_{L^p(bD)}\leq \|g^\#\|_{L^p(bD)}$. For $f\in L^1_{loc}(bD)$ and $g=f-\frac1{\vol(bD)} \int_{bD} f\d\sigma$, we have $g^\#=f^\#$ so 
\begin{align*}
\|M g\|_{L^p(bD)} \leq \|g^\#\|_{L^p(bD)}=\|f^\#\|_{L^p(bD)}. 
\end{align*}
 Since $Mf\leq Mg+\frac1{\vol(bD)}\left|\int_{bD} f\d\sigma\right|$,    
$\|Mf\|_{L^p(bD)}\leqs \|f^\#\|_{L^p(bD)}+\frac1{\vol(bD)} \left|\int_{bD} f\d\sigma\right|,$
which with (\ref{eq4}), proves Lemma \ref{Lemma2}.\qed
\par\medskip

From \cite{C-W}, $M$ is of weak-type $(1,1)$ and $(\infty,\infty)$. By Marcinkiewicz's theorem, $M:L^p(bD)\to L^p(bD)$ is continuous for all $1<p<+\infty$ which implies that $\|Mf\|_{L^p(bD)}\leqs \|f\|_{L^p(bD)}$ uniformly with respect to $f$. We also have $f^\#\leq 2Mf$, thus
\begin{align}
 \|f^\#\|_{L^p(bD)}&\leqs \|f\|_{\ch p},\quad 1<p<+\infty,\label{eq5}
\end{align}
and by definition
\begin{align}
 \|f^\#\|_{L^\infty(bD)}&\leqs \|f\|_{BMOA(D)}.\label{eq6}
\end{align}

Therefore, if $^\#$ was linear, its continuity (established in Lemma \ref{Lemma1}) would just be a simple consequence of interpolation.
Then, using Lemma \ref{Lemma2}, we would get 
$$\|f\|_{\ch q}\leqs\|f^\#\|_{L^q(bD)} \leqs \|f\|_{[\ch p ,BMOA(D)]_{1-\frac pq}}$$
as we need.
However, the operator $^\#$ is only sub-linear and it does not seem possible to get the continuity of $^\#:[\ch p,BMOA(D)]_{1-\frac pq} \to L^q(bD)$ that way. In \cite{FS}, Fefferman and Stein linearize the operator $^\#$. This technic requires a measurability condition that is not clearly satisfied in the case of convex domain of finite type because the extremal basis may have a chaotic behaviour. 
However, even if McNeal's polydiscs do not depend smoothly of $\zeta$ and $\varepsilon$, we will see that we  can have a smooth approximation of them with the Bergman metric because it is a smooth metric and because the ball centred at $\zeta$ of radius $1$ is almost equal to $\cp {|r(\zeta)|}\zeta$. This will allow us to define  another maximal function which will be comparable to $^\#$ and linearizable.

Without restriction, we assume that $0$ belongs to $D$ and we set
$$p(\zeta)=\inf\{\mu>0\ / \ \zeta \in\mu D\},$$ the gauge function of $D$, and $r=p-1$. 
We have $p(\lambda\zeta)=\lambda p(\zeta)$ for all $\lambda>0$. Therefore, for all $\varepsilon>0$, all $v\in\cc^n$, $|v|=1$, all $\lambda>0$,
 we have $\tau(\lambda \zeta,v,\lambda\varepsilon)=\lambda\tau(\zeta,v,\varepsilon)$.

We set, for $\zeta$ near $bD$ and $\varepsilon>0$ small, $\lambda(\zeta,\varepsilon):=\frac1{p(\zeta)-\varepsilon}$. When $\zeta$ is near $bD$, $p(\zeta)$ is near $1$ so $\lambda(\zeta,\varepsilon)$ is well defined. Moreover, we have $|r\big(\lambda(\zeta,\varepsilon)\zeta\big)|= \lambda(\zeta,\varepsilon)\varepsilon$, thus
$\tau\left(\lambda(\zeta,\varepsilon)\zeta,v,|r\big(\lambda(\zeta,\varepsilon)\zeta\big)|\right)=\lambda(\zeta,\varepsilon)\tau(\zeta,v,\varepsilon)$.

We now recall the definition and some properties of the Bergman metric that we will need (see \cite{Ran}). Let  $B(\zeta,z)$ denote the Bergman kernel,  holomorphic with respect to $z$, antiholomorphic with respect to $\zeta$ and let $b(z)=(b_{ij}(z))_{i,j=1\ldots, n}$ be  the matrix given by $b_{i,j}(z)=\diffp{^2}{z_i\partial\overline{z_j}} \ln B(z,z)$. The Bergman metric $\|\cdot\|_{B,z}$ for $z\in D$ is the hermitian metric induced by $b$, i.e. the Bergman norm of $v=\sum_{i=1}^n v_i e_i$, where $e_1,\ldots, e_n$ is the canonical basis of $\cc^n$, is given by $\|v\|_{B,z}=\left(\sum_{i,j=1}^n b_{ij}(z)v_j\overline{v_i}\right)^{\frac12}$.
\par
The two following  propositions are proved in \cite {McN} and \cite{McN2} respectively.
\begin{proposition}\mlabel{prop4}
 Let $\zeta\in D$ be a point near $bD$, $\varepsilon>0$, $w^*_1,\ldots, w^*_n$ a $\varepsilon$-extremal basis at $\zeta$ and $v= \sum_{j=1}^n v^*_j w^*_j$ a unit vector. Then, uniformly with respect to $\zeta, v$ and $\varepsilon $ we have
\begin{align*}
\frac1{\tau (z,v,\varepsilon )}&\eqs \sum_{j=1}^n\frac{|v^*_j|}{\tau_j(z,\varepsilon )}.
\end{align*}
\end{proposition}

\begin{proposition}\mlabel{prop2}
There exists $0<c<C$ such that for all $z\in D$ near $bD$, all unit vector $v$ in $\cc^n$
$$\frac c{\tau (z,v,|r(z)|)}\leqs \|v\|_{B,z}\leqs \frac C{\tau (z,v,|r(z)|)}.$$
\end{proposition}

Now we put for $\zeta$ near $bD$, $\varepsilon>0$
\begin{align*}
 \cq \varepsilon\zeta&=\left\{\zeta+\mu v\ / \ \mu\in\cc,\ v\in\cc^n,\ \|\mu v\|_{B,\lambda(\zeta,\varepsilon)\zeta}<1\right\},
\end{align*}
and, for $\alpha>0$ and $w_1^*,\ldots, w_n^*$ an $\varepsilon$-extremal basis at $\zeta$,
$$\alpha\cp\varepsilon\zeta=\left\{z=\zeta+\sum_{j=1}^nz_j^* w_j^*\ / \ |z_j^*|<\alpha\tau_j(\zeta,\varepsilon)\right\}.$$
Note that the factor $\alpha$ in front of $\cp\varepsilon\zeta$ means blowing up the polydisc around its center and not just multiplying each point by $\alpha$.
Now we prove 
\begin{proposition}\label{prop2.7}
 There exist $0<k<K$ such that for all $\zeta\in\cc^n$ near $bD$, all $\varepsilon>0$ small enough :
 $$k\cp\varepsilon \zeta \subset \cq \varepsilon \zeta \subset K \cp \varepsilon \zeta.$$
 \end{proposition}
\pr  We only prove the inclusion $k\cp\varepsilon \zeta \subset \cq \varepsilon \zeta $, the other one is similar. Let $v\in\cc^n$, $|v|=1$, $\mu\in\cc$ be such that $\zeta+\mu v$ belongs to $k\cp\varepsilon \zeta$, $k>0$ to be determined. From Proposition \ref{prop4}, we get
\begin{align*}
 \frac{|\mu|}{\tau(\zeta,v,\varepsilon)}\leqs k.
\end{align*}
Now $\tau(\lambda(\zeta,\varepsilon) \zeta,v,\lambda(\zeta,\varepsilon)\varepsilon)=\lambda(\zeta,\varepsilon) \tau(\zeta,v,\varepsilon)$  and 
$|r(\lambda(\zeta,\varepsilon)\zeta)|= \lambda(\zeta,\varepsilon)\varepsilon$
so
\begin{align*}
 \frac{|\mu|}{\tau(\lambda(\zeta,\varepsilon)\zeta,v,|r(\lambda(\zeta,\varepsilon)\zeta)|)}\leqs \frac k{\lambda(\zeta,\varepsilon)}.
\end{align*}
Therefore, by Proposition \ref{prop2}
$$\|\mu v\|_{B,\lambda(\zeta,\varepsilon) \zeta}\leqs \frac{k}{\lambda(\zeta,\varepsilon)}.$$
Since $\lambda(\zeta,\varepsilon)$ is near $1$, if $k$ is small enough, uniformly in $\zeta$ and $\varepsilon$, we have $\|\mu v\|_{B,\lambda(\zeta,\varepsilon)\cdot \zeta}<1$ and so $\zeta+\mu v$ belongs to $\cq \varepsilon \zeta$.\qed

\par\medskip 

Now we have all the tools we need in order to prove the continuity of $^\#$.

\begin{lemma}\mlabel{Lemma1}
For all $1<p<q<+\infty$, the operator  $^\#:[\ch p,BMOA(D)]_{1-\frac pq} \to L^q(bD)$ is continuous.
\end{lemma}
\pr  Let $f^*$ be defined by
$$f^*(z)=\sup_{\cq\varepsilon\zeta\ni z} \frac1{\vol(\cq\varepsilon \zeta \cap bD)}
\int_{\cq\varepsilon\zeta\cap bD}\left|f-f_{\cq\varepsilon\zeta\cap bD}\right|\d\sigma.$$
We show that the functions $f^*$ and $f^\#$ are comparable. Let $C>c>0$, depending only on the constants $k$ and $K$ given by Proposition \ref{prop2.7}, be such that for all $\zeta$ and $\varepsilon>0$, $\cp{c\varepsilon}\zeta\subset k\cp\varepsilon\zeta\subset \cq\varepsilon\zeta\subset K\cp\varepsilon \zeta\subset \cp{C\varepsilon}\zeta$. We have
\begin{align*}
 & \frac1{\vol(\cq\varepsilon \zeta \cap bD)}
\int_{\cq\varepsilon\zeta\cap bD}\left|f-f_{\cq\varepsilon\zeta\cap bD}\right|\d\sigma\\
&\leq
 \frac1{\vol(\cq\varepsilon \zeta \cap bD)}
\int_{\cq\varepsilon\zeta\cap bD}\left|f-f_{\cp{C\varepsilon}\zeta\cap bD}\right|\d\sigma+\left|f_{\cp{C\varepsilon}\zeta\cap bD}-f_{\cq\varepsilon\zeta\cap bD}\right|
\end{align*}
Since $\vol(\cq\varepsilon \zeta \cap bD)\eqs \vol(\cp{\varepsilon}\zeta\cap bD)\eqs \vol(\cp{C\varepsilon}\zeta\cap bD)$ and since 
$\int_{\cq\varepsilon\zeta\cap bD}|f-f_{\cp{C\varepsilon}\zeta\cap bD}|\d\sigma \leqs \int_{\cp{C\varepsilon}\zeta\cap bD}|f-f_{\cp{C\varepsilon}\zeta\cap bD}|\d\sigma$, we get
\begin{align}
\frac1{\vol(\cq\varepsilon \zeta \cap bD)}
\int_{\cq\varepsilon\zeta\cap bD}\left|f-f_{\cp{C\varepsilon}\zeta\cap bD}\right|\d\sigma &\leqs f^\#(z).\label{eq7} 
\end{align}
Now, using (\ref{eq7}), we get
\begin{align*}
 \left|f_{\cp{C\varepsilon}\zeta\cap bD}-f_{\cq\varepsilon\zeta\cap bD}\right|
 &=\left|\frac1{\vol(\cq\varepsilon \zeta \cap bD)}
\int_{\cq\varepsilon\zeta\cap bD}(f-f_{\cp{C\varepsilon}\zeta\cap bD})\d\sigma\right|\\
 &\leq  \frac1{\vol(\cq\varepsilon \zeta \cap bD)}
\int_{\cq\varepsilon\zeta\cap bD}\left|f-f_{\cp{C\varepsilon}\zeta\cap bD}\right|\d\sigma\\
&\leqs f^\#(z)
\end{align*}
which implies, again with (\ref{eq7}), that $f^*(z)\leqs f^\#(z)$. The converse inequality is analogue.

Now, since $\cq\varepsilon\zeta$ depends smoothly on $\zeta$ and $\varepsilon$, we can proceed as in \cite{FS} and linearize the maximal operator $^*$. Let $X$ be the set of couples $(\eta,Q)$ where $\eta:bD\times bD\to \cc$ is measurable and satisfies $|\eta|=1$ and 
$Q:bD\to\{ \cq\varepsilon \zeta\ /\ \varepsilon>0, \zeta\in bD\}$
is such that the map $(\zeta,z)\mapsto \mathds{1}_{Q(\zeta)}(z)$ is measurable and, for all $\zeta\in bD$, $\zeta$ belongs to $Q(\zeta)$. Let us point out that in order to define such a function $Q$, it suffices to define two functions $\phi:bD\to bD$ and $\psi:bD\to ]0,+\infty[$ and to set $Q(\zeta)=\cq{\psi(\zeta)}{\phi(\zeta)}$. Since $\mathds{1}_{Q(\zeta)}(z)=1$ if and only if $\left\|z-\phi(\zeta)\right\|_{B,\frac1{p(\phi(\zeta))-\psi(\zeta)} \phi(\zeta)}<1$ and since the Bergman metric is a smooth metric, $(\zeta,z)\mapsto \mathds{1}_{Q(\zeta)}(z)$ is  measurable as soon as that $\phi$ and $\psi$ are measurable. This would not be the case, a priori, with McNeal's polydiscs instead of $\cq\varepsilon\zeta$.

For $g\in L^1_{loc}(bD)$, we set $U(g)=\left(U_{(\eta,Q)}(g) \right)_{(\eta,Q)\in X}$ where
$$U_{(\eta,Q)}(g):\zeta \mapsto \frac1{\vol(Q(\zeta)\cap bD)}\int_{\xi\in Q(\zeta)\cap bD} \big(g(\xi)-g_{Q(\zeta)\cap bD}\big)\eta(\xi,\zeta) \d\sigma(\xi).$$
The operator $U$ is linear~; for all $(\eta,Q)\in X$, all $\zeta\in bD$, $\left|U_{(\eta,Q)}(g)(\zeta)\right|\leq g^*(\zeta)$ and for all $\zeta\in bD$, $\sup_{(\eta,Q)\in X}\left|U_{(\eta,Q)}(g)(\zeta)\right| = g^*(\zeta)$. In other words, for all $\zeta\in bD$, $U(g)(\zeta)$ is an element of $L^\infty(X)$ and $\|U(g)(\zeta)\|_{L^\infty(X)}=g^*(\zeta)$. We set 
$$L^\infty(bD,L^\infty(X))= \{f:bD\to L^\infty(X)\ /\ \|f\|_{L^\infty(bD,L^\infty(X))}=\sup_{bD}\|f\|_{L^\infty(X)}<\infty\},$$
$$L^p(bD,L^\infty(X))=\left\{f:bD\to L^\infty(X)\ /\ \|f\|_{L^p(bD,L^\infty(X))}:=\left(\int_{bD} \|f\|^p_{L^\infty(X)}\d\sigma \right)^{\frac1p}<\infty\right\}.$$
By definition $U:BMOA(bD)\to L^\infty(bD,L^\infty(X))$ is continuous. Moreover since $g^*\eqs g^\#$, from (\ref{eq5}) we deduce that $U:\ch p \to L^p(bD,L^\infty(D))$ is continuous. By interpolation, for all $\theta\in]0,1[$, 
$$U:[\ch p, BMOA(D)]_\theta\to \big[L^p(bD),L^\infty(X)),L^\infty(bD,L^\infty(X))\big]_\theta$$
is continuous. Since for $q$ such that $\frac1q=\frac{1-\theta}p+\frac\theta\infty$, $\big[L^p(bD,L^\infty(X)),L^\infty(bD,L^\infty(X))\big]_\theta=L^q(bD,L^\infty(X))$ (see Theorem 2.2.6 of \cite{HNV}), we conclude that uniformly with respect to $g\in [\ch p, BMOA(D)]_\theta$ : 
$$\left(\int_{bD} \|U(g)\|_{L^\infty(X)}^q\d\sigma\right)^{\frac1q} \leqs \|g\|_{[\ch p, BMOA(D)]_\theta}.$$
Since $\|U(g)\|_{L^\infty(X)}=g^*$, we get
$$\left(\int_{bD} {|g^*|}^q\d\sigma\right)^{\frac1q} \leqs \|g\|_{[\ch p, BMOA(D)]_\theta}$$
and finally, since $g^*\eqs g^\#$, this implies that $g^\#$ belongs to $L^q(bD)$ and satisfies, uniformly with respect to $g$,
$\|g^\#\|_{L^q(bD)}\leqs \|g\|_{[\ch p, BMOA(D)]_\theta}$. \qed
\par\medskip
We are now ready to prove the reciprocate of Lemma \ref{first_injection}.
\begin{lemma}\mlabel{Lemma3}
For $1<p< q<+\infty$, $[\ch p,BMOA(D)]_{1-\frac pq}\subset\ch q$ and for all $f\in  [\ch p,BMOA(D)]_{1-\frac pq}$, $\|f\|_{\ch q}\leqs \|f\|_{[\ch p,BMOA(D)]_{1-\frac pq}}$ uniformly with respect to $f$.
\end{lemma}
\pr 
For all $f\in[\ch p,BMOA(D)]_{1-\frac pq}$, we have (Lemma \ref{Lemma1}) :
\begin{align}
\|f^\#\|_{L^q(bD)}\leqs \|f\|_{[\ch p,BMOA(D)]_{1-\frac pq}}.\mlabel{eq2} 
\end{align}
In order to prove that $\frac 1{\vol(bD)}\left|\int_{bD} f\d\sigma\right|$ satisfies the same estimates, we consider the linear form $\lambda : L^1(bD)\to \cc$ defined by $\lambda(f)=\frac 1{\vol(bD)}\int_{bD} f\d\sigma$. The form $\lambda$ is continuous on $\ch 1$ and thus on $\ch p$ and $BMOA(D)$. Therefore, by interpolation, $\lambda$ is also continuous on $[\ch p,BMOA(D)]_{1-\frac pq}$ and for all $f\in[\ch p,BMOA(D)]_{1-\frac pq}$, we have 
\begin{align}
 |\lambda(f)|\leqs \|f\|_{[\ch p),BMOA(D)]_{1-\frac pq}}\mlabel{eq3}
\end{align}

Combining (\ref{eq2}) and (\ref{eq3}) with Lemma \ref{Lemma2}, we then get  for all $f\in[\ch p,BMOA(D)]_{1-\frac pq}$, $\|f\|_{L^q(bD)}\leqs \|f\|_{[\ch p,BMOA(D)]_{1-\frac pq}}$, so $[\ch p,BMOA(D)]_{1-\frac pq}$ injects itself continuously in $\ch q$.\qed 
\par\medskip 
Lemmas \ref{first_injection} and \ref{Lemma3} give immediately :
\begin{corollary}\mlabel{presque_interpolation}
For $1<p<q<+\infty$ and $\theta =1-\frac pq$, $[\ch p,BMOA(D)]_\theta=\ch q$ with equivalent norms.
\end{corollary}
\par\medskip
We now prove Theorem \ref{interpolation}.
\par\medskip
\noindent {\it Proof of Theorem  \ref{interpolation} :} 
First prove that $\ch q=[\ch 1,\ch p]_\theta$, for all $1<q<p<+\infty$ and $\theta$ such that $\frac1q=\frac{1-\theta}1+\frac\theta p$. Since $\ch p$ is reflexive, from  \cite{BL} Corollary 4.5.2, we have
$$[\ch 1,\ch p]'_\theta=[{\ch 1}',{\ch p}']_\theta.$$
We have ${\ch p}'=\ch{p'}$ where $\frac1p+\frac1{p'}=1$ and, from \cite{K-Li1},  $BMOA(D)={\ch 1}'$ so
\begin{align*}
[\ch 1,\ch p]'_\theta&=
 [BMOA(D),\ch {p'}]_\theta\\
 &=[\ch {p'},BMOA(D)]_{1-\theta}
\end{align*}
For $q'$ such that $\frac1{q'}=\frac\theta{p'}$, we have $1<p'<q'<+\infty$ and $1-\theta=1-\frac{p'}{q'}$. Thus Corollary \ref{presque_interpolation} implies that
\begin{align*}
[\ch 1,\ch p]'_\theta&=[\ch {p'},BMOA(D)]_{1-\theta}=\ch{q'}.
\end{align*}
Therefore, for $q$ such that $\frac1{q'}+\frac1q=1$ (which implies that $\frac1q=1-\theta+\frac\theta p$), we have

\begin{align*}
[\ch 1,\ch p]''_\theta&=\ch{q}.
\end{align*}

Since $[\ch 1,\ch p]_\theta$ is a subspace of $[\ch 1,\ch p]''_\theta$ which is reflexive since isomorphic to the reflexive space $\ch q$, it follows that $[\ch 1,\ch p]_\theta$ is itself reflexive and so $[\ch 1,\ch p]_\theta=\ch{q}$ and $\frac1q=\frac{1-\theta}1+\frac\theta p$.

Now we prove that $\ch q=[\ch 1,BMOA(D)]_{1-\frac1q}$. For $1<q<+\infty$, $\theta=\frac12$, $\theta'=\frac{2(q-1)}{2q-1}=1-\frac1{2q-1}\in ]0,1[$, we have $\frac1q=1-\theta'+\frac{\theta'}{2q}$ and so
\begin{align*}
 [\ch q,BMOA(D)]_\theta&= \ch {2q},\\
 [\ch 1,\ch {2q}]_{\theta'}&= \ch q.
\end{align*}
Therefore, for $s=\frac{\theta\theta'}{1-\theta'+\theta\theta'}$, we get from Wolff's note \cite{Wol}, Theorem 2 :
$$[\ch 1,BMOA(D)]_s= \ch q.$$
Since $s=1-\frac 1q$, we are done.\qed

\section{About $h\omega$ and $w$}\mlabel{bazar}
\subsection{$h\omega$ satisfies the hypothesis of Theorem \ref{ddbar}}\mlabel{subsection31}
In this section, we prove that $h\omega=\frac h{(|f_1|^2+|f_2|^2)^2}\left(\overline{f_1}\,\overline\partial \, \overline{f_2} - \overline{f_2}\,\overline \partial\, \overline{f_1}\right)$, as defined in the introduction, satisfies the hypothesis of Theorem \ref{ddbar}. 
When $\theta$ is a smooth $p$-form on $\overline D$, Bruna, Charpentier and Dupain \cite{BCD} define $\|\theta(z)\|_{k}$ as a smooth function of $z$ by $\|\theta(z)\|_{k}:=\sup_{u_1,\ldots, u_p\neq 0} \frac{|\omega(z)(u_1,\ldots, u_p)|}{k(z,u_1)\ldots k(z,u_p)}$ which is the norm of the form $\theta(z)$ with respect to the norm $k(z,\cdot)$. The following theorem
is Theorem 1.2 of \cite{Jas} except for the estimates 
$\left\||r|\, \|\partial h\wedge \overline{\partial h}\|_k\d V\right\|_{W^1}\leqs \|h\|^2_{BMOA(D)}$ and $\left\||r|\, \|\partial h\|_k^2\d V\right\|_{W^1}\leqs \|h\|^2_{BMOA(D)}$ which are not stated in it, but they are immediate consequences of the inequality $\int_{\cp \varepsilon z \cap bD}|r|\, \left\|\partial h \wedge \overline{\partial h}\right\|_k \leqs \|h\|_{BMOA(D)}^2\sigma(\cp \varepsilon z\cap bD)$ established for all $\varepsilon>0$ and all $z\in bD$ in the proof of Theorem 1.2 of \cite{Jas}.

\begin{theorem}\mlabel{fact2lemma2}
For all $h\in BMOA(D)$, $|r|\ \|\partial h\wedge \overline{\partial h}\|\d V$ and $|r|\, \|\partial h\|_k^2\d V$ are Carleson measure and  
$\big\||r|\, \|\partial h\wedge \overline{\partial h}\|_k\d V\big\|_{W^1}\leqs\hskip -1pt \|h\|_{BMOA(D)}^2$ and $\big\||r|\,\|\partial h\|_k^2\d V\big\|_{W^1}\leqs \hskip -1pt \|h\|_{BMOA(D)}^2$.
\end{theorem}
The following  theorem is Theorem 1.1 of \cite{Jas}. This result, generally referred as Carleson-Hörmander inequality, will be very useful for us. 

\begin{theorem}[Carleson-Hörmander inequality] \mlabel{Carleson-Hormander}
 Let $D$ be a bounded convex domain of finite type in $\cc^n$, let $\mu$ be a Carleson measure in $D$. Then for all $1<p<+\infty$, all $h\in\ch p$, we have, uniformly with respect to $h$ :
 $$\int_D|h|^p\d\mu\leq\int_{bD} |h|^p\d\sigma.$$
\end{theorem}

We will also need the following lemma which is the analog of Proposition 2.1 of \cite{AndeCarl1}.
\begin{lemma}\mlabel{fact2lemma1}
 For all $p\in]0,+\infty[$, all $h\in \ch p$, we have
 \begin{align*}
  \int_D |r|\, |h|^{p-2} \|\partial h\|_k^2 \d V&\leqs \|h\|^p_{\ch p}.
 \end{align*}
\end{lemma}
\pr
We put $\theta=i\partial\overline\partial |h|^p=i\left(\frac{p}2\right)^2 |h|^{p-2}\partial h\wedge \overline{\partial h}$.

Since $i\partial h\wedge\overline{\partial h}(v,iv)=2|\partial h(v)|^2$, we have
\begin{align*}
 |h|^{p-2}\|\partial h\|_k^2 &\leqs |h|^{p-2}\|\partial h \wedge \overline{\partial h}\|_k \eqs\|\theta\|_k,
\end{align*}
so
\begin{align*}
 \int_D |r|\, \|h\|^{p-2}\|\partial h\|_k^2\d V\leqs \int_D |r|\, \|\theta\|_k\d V.
\end{align*}
Since $\theta$ is a  closed positive $(1,1)$-current, Theorem 1.1 of \cite{BCD} gives
\begin{align*}
 \int_D |r|\, \|h\|^{p-2}\|\partial h\|_k^2\d V&\leqs \int_D |r|\, \|\theta\|_{eucl}\d V.
\end{align*}
where $\|\theta\|_{eucl}$ stands for the euclidean norm of $\theta$ :
\begin{align*}
 \|\theta(z)\|_{eucl}=\sup_{u,v\neq 0} \frac{|\theta(z)(u,v)|}{|u|\cdot |v|}.
\end{align*}
Let $e_1,\ldots, e_n$ be the canonical basis of $\cc^n$ and let us write $u=\sum_{j=1}^n u_je_j$ and $v=\sum_{j=1}^n v_je_j$. Since $\theta$ is positive
\begin{align*}
 |\theta(z)(u,v)|\leq \sqrt{\theta(z)(u,iu)} \sqrt{\theta(z)(v,iv)}
\end{align*}
which yields
\begin{align*}
 \frac{|\theta(z)(u,v)|}{|u| \cdot|v|} &\leq \sum_{k,j=1}^n \theta(z)(e_j,e_k)\frac{u_jv_k}{|u|\, |v|}\\
 &\leqs \sum_{j=1}^n \theta(e_j,ie_j)\\
 &\leqs \Delta |h|^p
\end{align*}
so
\begin{align*}
 \int_D |r|\, \|h\|^{p-2}\|\partial h\|_k^2\d V\leqs \int_D |r|\, \Delta |h|^p\d V.
\end{align*}
Now by Green identity
\begin{align*}
 \int_D|r|\, \Delta |h|^p\d V +\int_D |h|^p \Delta r \, \d V= \int_{bD} |h|^p \diffp r \eta \d\sigma
\end{align*}
and since $r$ is convex, $\Delta r\geq 0$ so
\begin{align*}
 \int_D |r|\, \|h\|^{p-2}\|\partial h\|_k^2\d V &\leqs \int_{bD} |h|^p \diffp r \eta \d\sigma\\
 &\eqs \int_{bD} |h|^p\d\sigma=\|h\|_{\ch p}^p.\qed
\end{align*}

Since $\overline\partial(h\omega)=0$, $\partial(h\omega)$ is $\d$-closed. Thus Fact \ref{fact2} below shows that $\partial( h\omega)$ satisfies the hypothesis of Theorem \ref{ddbar}. We first establish :

\begin{fact}\mlabel{fact1}
For $f_1, f_2\in \ch\infty$ such that $|f_1|^2+|f_2|^2\geq\delta^2>0$ and $\omega=\frac{\overline{f_1}\,\overline\partial \, \overline{f_2} - \overline{f_2}\,\overline \partial\, \overline{f_1}} {(|f_1|^2+|f_2|^2)^2}$, 
$|r|\, \|\partial \omega\|_k\d V$ and  $|r|\, \|\omega\|_k^2\d V$ are Carleson measures on $D$.
\end{fact}
\pr 
We have $\|\omega\|_k\leqs \|\partial f_1\|_k+\|\partial f_2\|_k$ so
\begin{align}
 |r|\, \|\omega\|_k^2\leqs |r|\, (\|\partial f_1\|^2_k+\|\partial f_2\|^2_k).\mlabel{fact1eq1}
\end{align}
We also have
\begin{align*}
 \partial \omega&= \frac2{(|f_1|^2+|f_2|^2)^3} \left( \overline{f_1}^2 \overline{\partial f_2}\wedge {\partial f_1}-\overline{f_2}^2 \overline{\partial f_1}\wedge {\partial f_2} +\overline{f_1f_2}\left( \overline{\partial f_2}\wedge \partial f_2-\overline{\partial f_1}\wedge \partial f_1\right)\right)
\end{align*}
and since for all $\alpha$ and $\beta$, $\|\alpha\wedge\beta\|_k\leqs\|\alpha\|_k\|\beta\|_k$, we get
\begin{align}
 |r|\,\|\partial \omega\|_k&\leqs |r|(\|\partial f_1\|_k^2+\|\partial f_2\|_k^2).\mlabel{fact1eq2}
\end{align}
Since $\ch\infty\subset BMOA(D)$, Fact \ref{fact1} is then a consequence of (\ref{fact1eq1}), (\ref{fact1eq2}) and of Theorem \ref{fact2lemma2}.\qed

\begin{fact}\mlabel{fact2}
 For all $p\in [1,+\infty]$, all $h\in\ch p$, $|r|\partial (h\omega)$ and $|r|\partial h \wedge \omega $ belong to $W^{1-\frac1p}_{1,1}(D)$. 

\end{fact}
\pr 
We treat separately the case $h\in \ch \infty$. Since $\partial (h\omega)=\partial h \wedge \omega+h\partial \omega$, it suffices to prove that both $|r| \partial h \wedge \omega$ and $ |r| h\partial \omega$ belongs to $W^1_{1,1}(D)$. Since for all vectors fields $u_1$ and $u_2$ we have  $\frac1{k(\cdot, u_1)k(\cdot,u_2)}|\partial h\wedge \omega (u_1,u_2)|\leq \|\partial h\wedge \omega \|_k$ and $\frac1{k(\cdot, u_1)k(\cdot,u_2)}|\partial \omega (u_1,u_2)|\leq \|\partial \omega \|_k$, we just prove that $|r|\, \|\partial h\wedge\omega\|_{k}\d V$ and  $|r|\, \|h\partial\omega\|_{k}\d V$ are Carleson measure. 

From fact \ref{fact1}, $|rh|\, \|\partial \omega \|_k\d V$ is a Carleson measure.
Next, since $\|\partial h\wedge \omega\|_k\leqs\|\partial h\|_k\|\omega\|_k$  for any $\zeta_0\in bD$ and all $\varepsilon>0$, we get from Theorem \ref{fact2lemma2} and Fact \ref{fact1} :
\begin{align*}
 \int_{\cp \varepsilon{\zeta_0}\cap D} |r|\ \|\partial h\wedge \omega\|_k\d V
 &\leqs\left(\int_{\cp \varepsilon{\zeta_0}\cap D} |r|\ \|\partial h\|_k^2\d V \right)^{\frac12}
 \left(\int_{\cp \varepsilon{\zeta_0}\cap D} |r|\ \|\omega \|_k^2\d V \right)^{\frac12}\\
 &\leqs \|h\|_{BMOA(D)}\sigma\big(\cp \varepsilon{\zeta_0}\cap D \big)\\
 &\leqs \|h\|_{\ch\infty}\sigma\big(\cp \varepsilon{\zeta_0}\cap D \big).
\end{align*}

Now we treat the case $h\in \ch p$, $p\in[1,+\infty[$. 
It suffices to prove that  $|r|\,\partial h \wedge \omega$ belongs to $W_{1,1}^{1-\frac1p}$ and $|r|\,|h|\|\partial \omega \|_k$ belong to $W^{1-\frac1p}(D)$.

From Fact \ref{fact1}, $|r|\, \|\partial \omega\|_k\d V$ is a Carleson measure. The Hörmander-Carleson inequality implies that $h$ belongs to $L^p(D, |r|\, \|\partial \omega\|_k\d V)$. Then, Proposition \ref{propAB}, $|rh|\, \|\partial \omega\|_k\d V$ belongs to $W^{1-\frac1p}(D)$.

Now it remains to prove that  $|r|\,\partial h \wedge \omega$ belongs to $W_{1,1}^{1-\frac1p}(D)$. We proceed by interpolation. Let us consider the linear operator $T:h\mapsto \partial h\wedge \omega$. We prove that $T:\ch1\to W_{1,1}^0(D)$ and $T:BMOA(D)\to W^1_{1,1}(D)$ are continuous.

Let $h$ belongs to $\ch1$ and let $u_1$ and $u_2$ be two vectors fields. Then :
\begin{align*}
 \int_{D}\frac {|r|\,\left|\partial h\wedge \omega (u_1,u_2)\right|} {k(\cdot,u_1)k(\cdot,u_2)} 
 &\leqs \int_D |r|\ \|\partial h\wedge \omega\|_k\d V.
\end{align*}
Since $\|\partial h \wedge \omega \|_k\leqs \|\partial h\|_k \cdot\|\omega\|_k$, we have
\begin{align*}
 \int_D|r|\, \|\partial h\wedge \omega\|_k\d V
 &\leq \left(\int_D  |r| |h|^{-1} \|\partial h\|_k^2 \d V\right)^{\frac12}
  \left(\int_D |r| {|h|}\|\omega\|_k^2\d V\right)^{\frac12}.
\end{align*}
If $h$ belongs to $\ch 1,$ Lemma \ref{fact2lemma1} implies that $\int_D  |r|\, |h|^{-1} \|\partial h\|_k^2 \d V \leqs\|h\|_{\ch 1}.$

Since, Fact \ref{fact1}, $|r|\,\|\omega\|_k^2$ is a Carleson measure, if $h$ belongs to $\ch 1$, Hörmander-Carleson inequality yields $\int_D |r| {|h|}\|\omega\|_k^2\d V\leqs \|h\|_{\ch1}$ and so
\begin{align*}
 \int_D|r|\, \|\partial h \wedge \omega\|_k^2\d V\leqs \|h\|_{\ch 1}.
\end{align*}
Thus, we have proved that $T:\ch 1 \to W_{1,1}^0(D)$ is continuous.

Now for $h\in BMOA(D)$, $u_1$ and $u_2$ be two vectors fields, $z_0\in bD$, $\varepsilon>0$, we have
\begin{align*}
 \int_{\cp\varepsilon {z_0}\cap bD}\frac{ |r| \left|\partial h\wedge \omega (u_1,u_2)\right| }{k(\cdot,u_1)k(\cdot,u_2)}
 &\leqs \int_{\cp\varepsilon {z_0}\cap bD} |r|\,\|\partial h\wedge \omega\|_k\d V\\
 &\leqs \left(\int_{\cp\varepsilon {z_0}\cap bD} |r|\,\|\partial h\|_k^2\d V \right)^{\frac12}
 \left(\int_{\cp\varepsilon {z_0}\cap bD} |r|\,\|\omega\|_k^2\d V\right)^{\frac12}
\end{align*}
From Theorem \ref{fact2lemma2} and Fact \ref{fact1}, we then obtain :
\begin{align*}
 \int_{\cp\varepsilon {z_0}\cap bD}|r|\,\|\partial h \wedge \omega\|_k\d V &\leqs \|h\|_{BMOA(D)} \sigma({\cp\varepsilon {z_0}\cap bD}).
\end{align*}
We thus have proved that $T:BMOA(D)\to W^1_{1,1}(D)$ is continuous.

By interpolation, we get the continuity of $T: \ch p\to W^{1-\frac1p}_{1,1}(D)$ for all $1\leq p<\infty$. Thus, for all $h\in\ch p$, $1\leq p<+\infty$, $|r|\, \partial h\wedge \omega $ belongs to $W^{1-\frac1p}_{1,1}(D)$.\qed

\subsection{$w$ belongs to $\ch p$}\mlabel{subsection32}
Using that $w$ is a holomorphic function such that $\overline w(0)=v(0)$ and $\d\overline w=\overline\partial v-h\omega$ with $v\in L^p(bD)$ and $h\omega\in \ch p$, we now prove that $w$ is in $\ch p$. Let $p'$ be such that $\frac1p+\frac1{p'}=1$. We test $w$ against any function $g\in \ch {p'}$ and showing that $\left|\int_{bD} g\overline w d\sigma \right| \leqs \|g\|_{\ch {p'}}$, uniformly with respect to $g$, we get by duality that $w$ belongs to $\ch p$. Moreover, since $v$ belongs to $L^p(bD)$, it suffices to prove that  $\left|\int_{bD} g(\overline w-v) d\sigma \right| \leqs \|g\|_{L^{p'}(bD)}$ uniformly with respect to $g$. We use the following lemma.
\begin{lemma}\label{difficile}
For $1<p\leq+\infty$ and $1\leq p'<+\infty$ such that $\frac1p+\frac1{p'}=1$, $g\in\ch{p'}$, $h\in\ch p$, we have uniformly with respect to $g$ and $h$ :
\begin{align*}
 \int_D|r|\,|\partial g|\, |h|\, |\omega|\d V&\leqs \|h\|_{\ch p}\|g\|_{\ch{p'}}.
\end{align*}
\end{lemma}
\pr The proof is similar to a part of the proof of Theorem 1.2 of \cite{AndeCarl1}.

If $p=+\infty$ and $p'=1$, 
from Fact \ref{fact1}, $|r|\, \|\omega\|_k^2\d V$ is a Carleson measure. Thus 
from Lemma \ref{fact2lemma1} and Theorem \ref{Carleson-Hormander} applyed to $g$ and $\mu=|r|\,\|\omega\|_k^2\d V,$ we get :
\begin{align*}
 \int_D|r|\,|\partial g| |h| |\omega|\d V&
 \leqs \|h\|_{\ch\infty}\left(\int_D |r|\ |g|^{-1}\|\partial g\|_k^2\d V\right)^{\frac12}\left(\int_D |r| |g| \|\omega\|_k^2\d V\right)^{\frac12}\\ 
 &\leqs \|h\|_{\ch \infty}\|g\|_{\ch{1}}.
\end{align*}
If $p=p'=2,$ again from Lemma \ref{fact2lemma1}, Fact \ref{fact1} and Carleson-Hörmander inequality applyed to $h$ and $\mu=|r|\,\|\omega\|_k^2\d V,$, we have :
\begin{align*}
 \int_D|r|\,|\partial g| |h| |\omega|\d V&
 \leqs\left(\int_D |r|\ \|\partial g\|_k^2\d V\right)^{\frac12}\left(\int_D |r| |h|^2 \|\omega\|_k^2\d V\right)^{\frac12}\\ 
 &\leqs \|h\|_{\ch 2}\|g\|_{\ch{2}}.
\end{align*}
If $2<p<+\infty$, $1<p'<2$, since $\frac{p'}{2-p'}$ and $\frac p2$ are dual exponents, still from Lemma \ref{fact2lemma1}, Fact \ref{fact1} and Carleson-Hörmander inequality :
\begin{align*}
 &\int_D|r|\,|\partial g| |h| |\omega|\d V\\
 &\leqs\left(\int_D |r|\ |g|^{p'-2} \|\partial g\|_k^2\d V\right)^{\frac12}
 \left(\int_D |r| |g|^{2-p'}|h|^2\|\omega\|_k^2\d V\right)^{\frac1{2}}\\
 &\leqs\left(\int_D |r|\ |g|^{p'-2} \|\partial g\|_k^2\d V\right)^{\frac12}
 \hskip -6pt\left(\int_D |r| |g|^{p'}\|\omega\|_k^2\d V\right)^{\frac{2-p'}{2p'} }
 \hskip -6pt \left(\int_D |r| |h|^{p}\|\omega\|_k^2\d V\right)^{\frac1{p}} \\ 
 &\leqs \|h\|_{\ch p}\|g\|_{\ch{p'}}.
\end{align*}
Finally, if $1<p<2,$ $2<p'<+\infty$, we write $h\partial g$ as $\partial(hg)-g\partial g$. Applying the case $p=+\infty$, $p'=1$ to the function identically equal to  $1$ on $D$ and to $gh\in\ch 1$, we get 
\begin{align}
\int_D|r|\,|\partial (gh)| |\omega|\d V&\leqs \|hg\|_{\ch1}\leq\|g\|_{\ch {p'}} \|h\|_{\ch {p}}.\label{eqcor1}
\end{align}
Applying the case $2<p<+\infty$, $1<p'<2$ and to $h\in\ch{p}$ and $g\in\ch {p'}$, we get 
\begin{align}
\int_D|r|\,|\partial h| |g| |\omega|\d V\leqs \|g\|_{\ch {p'}} \|h\|_{\ch {p}}. \label{eqcor2}
\end{align}
Inequalities (\ref{eqcor1}) and (\ref{eqcor2}) imply that  $\int_D|r|\,|\partial g| |h| |\omega|\d V\leqs \|h\|_{\ch p}\|g\|_{\ch{p'}}$ when $1<p<2,$ $2<p'<+\infty$.\qed

Now let $G$ be the Green function for the Laplacian for $D$ (see \cite{Kra}). We will need the following properties of $G$ :
\begin{itemize}
 \item The following representation formula holds for all $f\in C^2(\overline D)$ and all $z\in D$, 
 $$f(z)=-\int_{bD}f(\xi)\diffp G{\eta_\xi}(z,\xi)\d\sigma(\xi)
 +\int_DG(z,\xi)\Delta f(\xi)\d V(\xi).$$
\item For all $z\in D$, all $\zeta\in \overline{D},$ $z\neq\zeta$, $G(z,\zeta)\geq 0$. Indeed, let us consider $\varphi(\zeta)=-G(z,\zeta)$, $z$ fixed in $D$. For all $\varepsilon>0$, $\varphi$ is harmonic in $D\setminus B(z,\varepsilon)$, $\varphi(\zeta)=0$ on $bD$ and $\varphi(\zeta)<0$ for all $\zeta\in bB(z,\varepsilon)$, all $\varepsilon>0$ sufficiently small. The maximum principle implies that $\varphi(\zeta)\leq 0$ for all $\zeta\in D\setminus B(z,\varepsilon)$, $\varepsilon>0$ sufficiently small, so $G(z,\zeta)\geq 0$ for all $z\in D$, all $\zeta\in\overline D$, $\zeta\neq z$.
\item For all $z\in D$, all $\zeta\in bD$, $\diffp G{\eta_\zeta}(z,\zeta)<0$. Indeed, for all $z\in D$ fixed, $\varphi:\xi\mapsto G(z,\xi)$ is harmonic in $D\setminus B(z,\varepsilon)$, $\varepsilon>0$. For all $\xi\in D\setminus\{z\}$, all $\zeta\in bD$ fixed, we have $\varphi(\xi)\geq \varphi(\zeta)=0$. Hopf's lemma (see \cite{Krabis}) implies that $\diffp \varphi {\eta_\zeta}(\zeta)<0$.
 \end{itemize}
We put $G_0=G(0,\cdot)$. Since $-\diffp{G_0}{\eta}>0$ on the compact set $bD$, we have 
\begin{align*}
 \left|\int_{bD}g(\overline w-v)\d\sigma\right|&\eqs \left|\int_{bD}g(\overline w-v)\diffp{G_0}{\eta} \d\sigma\right|. 
\end{align*}
Since $\overline w(0)=v(0)$, the representation formula gives 
\begin{align*}
 \left|\int_{bD}g(\overline w-v)\d\sigma\right|&\eqs \left|\int_{D}G_0 \Delta( g(\overline w-v)) \d V\right|. 
\end{align*}
Now let $\beta=\frac i2 \sum_{k=1}^n\d z_k\wedge\d \overline{z_k}$. For any $f$ we have $\partial\overline\partial f\wedge\beta^{n-1}\eqs c_n \Delta f \d V$ for some $c_n\in\cc$ depending only on $n$. So, since $g$ and $w$ are holomorphic, since $\partial\overline\partial v=\partial (h\omega)$ and $\overline\partial (\overline w -v)=-h\omega$ :
\begin{align*}
 \left|\int_{bD}g(\overline w-v)\d\sigma\right|&\eqs \left|\int_{D}G_0 \partial\overline\partial( g(\overline w-v)) \wedge \beta^{n-1}\right|\\
 &\leqs  \left|\int_{D}G_0 \partial g\wedge\overline{\partial}(\overline w-v) \wedge \beta^{n-1}\right|
 +\left|\int_{D}G_0  g \partial\overline{\partial}v \wedge \beta^{n-1}\right|
 \\
 &\leqs \left|\int_{D}G_0 h\partial g\wedge \omega \wedge \beta^{n-1}\right|
 +\left|\int_{D}G_0  g \partial(h\omega) \wedge \beta^{n-1}\right|.
\end{align*}
Let $K$ be a compact neighborhood of $0$. Since $G_0$ is of class $C^1(\overline D\setminus K)$ and vanishes on $bD$, $|r|^{-1} G_0$ is bounded on $D\setminus K$. Thus, by Lemma \ref{difficile} :
\begin{align*}
 \left|\int_{D\setminus K} G_0h\partial g\wedge \omega \wedge \beta^{n-1}\right|
 &\leq \int_{D\setminus K} \frac{G_0}{|r|} |h| | \partial g|\, |r|\,|\omega| \d V\\
 &\leqs \|g\|_{\ch {p'}} \|h\|_{\ch p}.
\end{align*}
Since $G_0$ is locally integrable and since $\omega$ is bounded on $K$ :
\begin{align*}
 \left|\int_KG_0h\partial g\wedge\omega\wedge\beta^{n-1}\right|
 &\leqs \sup_K|\partial g|\sup_K|h|\\
 &\leqs \|g\|_{L^{p'}(D)}\|h\|_{L^p(D)}\\
 &\leqs \|g\|_{\ch{p'}}\|h\|_{\ch p}
\end{align*}
and so $\left|\int_DG_0h\partial g\wedge\omega\wedge\beta^{n-1}\right|\leqs \|g\|_{\ch{p'}}\|h\|_{\ch p}$.

Now we show that $\left|\int_{D}G_0  g \partial(h\omega) \wedge \beta^{n-1}\right|\leqs\|g\|_{\ch p}$. 
Let $e_1,\ldots, e_n$ be the canonical basis of $\cc^n$. We have
\begin{align*}
 \left|\int_{D\setminus K}G_0  g \partial(h\omega) \wedge \beta^{n-1}\right|&\leqs \sum_{i,j=1}^n \left|\int_{D\setminus K} \frac {G_0}{|r|} g \frac{|r| \partial(h\omega)(e_i,e_j)}{k(\cdot,e_i)k(\cdot,e_j)} \d V\right|.
\end{align*}
From Fact \ref{fact2}, $|r| \frac{\partial(h\omega)(e_i,e_j)}{k(\cdot,e_i)k(\cdot,e_j)}$ belongs to $W^{1-\frac1p}(D)$ so, Proposition \ref{propAB},  it can be written as $f\d\mu$ where $\mu$ is a Carleson measure on $D$ and $f$ belongs to $L^p(D,\mu)$. Therefore, since $\frac{G_0} {|r|}$ is bounded on $D\setminus K$ :
\begin{align*}
 \left|\int_{D\setminus K} \frac {G_0}{|r|} g \frac{|r| \partial(h\omega)(e_i,e_j)}{k(\cdot,e_i)k(\cdot,e_j)} \d V\right|
 &\leqs \left|\int_{D\setminus K}  g f \d \mu\right|\\
 &\leqs \left(\int_{D\setminus K}  |g|^{p'}\d\mu\right)^{\frac1{p'}}  \left(\int_{D\setminus K} |f|^p \d \mu\right)^{\frac1p}.
 \end{align*}
 and using Carleson-Hörmander inequality, we then get :
 \begin{align*}
 \left|\int_{D\setminus K} \frac {G_0}{|r|}g \frac{{|r|}\partial(h\omega)(e_i,e_j)}{k(\cdot,e_i)k(\cdot,e_j)} \d V\right|
 &\leqs\|g\|_{\ch {p'}}.
\end{align*}
and so $ \left|\int_{D}G_0  g \partial(h\omega) \wedge \beta^{n-1}\right|\leqs \|g\|_{\ch{p'}}$.

Since $G_0$ is locally integrable, as previously, we have :
\begin{align*}
 \left|\int_{K}G_0  g \partial(h\omega) \wedge \beta^{n-1}\right|&\leqs \sup_{K}|g| \sup_{K}|h\omega|\\
 &\leqs \|g\|_{\ch{p'}} \|h\|_{\ch p}.
\end{align*}
So, for any $g\in\ch{p'}$, $\left|\int_{bD} g\overline w d\sigma\right|\leqs \|g\|_{\ch {p'}}$ which implies that $w$ belongs to $\ch p$.

\section{Proof of Theorem \ref{ddbar}}\mlabel{section5}
 The proof of Theorem \ref{ddbar} reduces to the 2 following theorems.
\begin{theorem}\mlabel{poincare}
Let $D$ be a bounded convex domain with smooth boundary of finite type, let  $\theta$ be a closed positive $(1,1)$-current such that $|r|\theta$ belongs to $W^\alpha_{1,1}(D)$ for some $\alpha\in[0,1]$.\\
Then there exists $v$ real $1$-form in $W^\alpha_1(D)$ such that $\d v=\theta$ and 
$$\|v\|_{W^\alpha_1}(D)\leqs \||r|\, \theta\|_{W^\alpha_{1,1}(D)},$$
uniformly with respect to $\theta$.
\end{theorem}
\begin{theorem}\mlabel{dbar}
 Let $D$ be a bounded convex domain with smooth boundary of finite type. For all $\overline\partial$-closed $v\in C^\infty_{0,1}(\overline D)\cap W_{0,1}^{1-\frac1p}(D)$, $p\in[1,+\infty]$, there exists $u\in C^\infty(\overline D)$ such that
 \begin{itemize}
  \item $\overline\partial u=v$,
  \item $\|u\|_{L^p(bD)}\leqs \|v\|_{W_{0,1}^{1-\frac1p}(D)}$ if $1\leq p<+\infty$,
  \item $\|u\|_{BMO(bD)} \leqs \|v\|_{W^1_{0,1}(D)}$ if $p=+\infty$.
 \end{itemize}
\end{theorem}
Theorem \ref{dbar} is Theorem 2.10 of \cite{Ale2}. Theorem \ref{poincare} will be proved by interpolation. We admit it for the moment and prove Theorem \ref{ddbar}.
\par\medskip
{\it Proof of Theorem \ref{ddbar} :} This is classic, we include it for completness. Since $\theta$ is positive, it is real and since $\theta$ is $\d$-closed, there exists $v$ real $1$-form such that $i\d v=\theta$. We decompose $v=-v_{1,0}+v_{0,1}$ where $v_{0,1}$ is a $(0,1)$-form and $v_{1,0}$ a $(1,0)$-form. For bidegree reason $\overline\partial v_{0,1}=0$. Let $u$ be such that $\overline\partial u=v_{0,1}$. We put $w=2\Re u$ and, using $\overline v_{0,1} =v_{1,0}$, we get
\begin{align*}
 i\partial\overline\partial w
 &=i\partial\overline\partial (u+\overline u)\\
 &=i\partial\overline\partial u-i\overline\partial\partial\overline u\\
 &=i\partial v_{0,1} -i\overline\partial \overline{v_{0,1}}\\
 &=i\d v=\theta.
 \end{align*}
Now when $v$ is given by Theorem \ref{poincare}, $v_{0,1}$ belongs to $W^{1-\frac1p}_{0,1}(D)$ if $|r|\,\theta$ belongs to $W^{1-\frac1p}_{1,1}(D)$ and then, when $u$ is given by Theorem \ref{dbar}, $w$ belongs to $BMO(bD)$ if $p=+\infty$ and to $L^p(bD)$ if $1\leq p< +\infty$.\qed

Our goal is now to prove Theorem \ref{poincare}. 
We will use the homotopy operator of \cite{Ale3} that we now recall. Let $\varphi $ be a $C^\infty$ smooth function such that $\varphi (t)=1$ if $t<\frac12$, $\varphi (t)=0$ if $t>1$, and define the map $h_\Lambda :D\times [0,1]\to D$ for $|\Lambda|\leq \rho$ by
\begin{align*}
h_\Lambda (z,t)
&=tz+ t\varphi \left(\frac{1-t}{\gamma|r(z)|}\right)\frac{1-t}{|r(z)|} A(z) \cdot\Lambda   +t\left(1-\varphi \left(\frac{1-t}{\gamma|r(z)|}\right)\right)A(tz) \cdot\Lambda 
\end{align*}
where $\gamma$ and $\rho$ have to be chosen sufficiently small, $A(z)$ is a positive hermitian matrix, smoothly depending on $z$, such that $A(z)^{-2}=B(z)$, $B(z)$ being the matrix in the canonical basis which determines the Bergman metric $\|\cdot\|_{B,z}$ at $z$, i.e. $\|v\|_{B,z}=\overline v^t B(z) v$ for any vector $v$. The map $h_\Lambda $ is $C^\infty$-smooth in $D\times ]0,1[$, $h_\Lambda (z,0)=0$ and $h_\Lambda (z,1)=z$ for all $z$ in $D$.\\
The associated homotopy operator is
$$H\theta  =\frac1{\vol(\ \Delta_n(\rho))}\int_{\Lambda \in \Delta_n(\rho)} \left(\int_{t\in [0,1]} h^*_\Lambda \theta  \right)\d\Lambda,$$
where $\Delta_n(\rho)=\{\Lambda\in\cc^n,\ |\Lambda|<\rho\}$.
If $\theta$ is closed and if its support does not meet $0$, then $\mathrm{d}H\theta=\theta$.

Moreover, the author proved in \cite{Ale3} that for all closed positive $(1,1)$-current $\theta$ supported away from the origin and such that $|r|\,\theta$ belongs to $W^1_{1,1}(D)$, $H\theta$ belongs to $W^1_{1}(D)$ and satisfies $\|H\theta\|_{W^1_1(D)}\leqs \||r|\, \theta\|_{W^1_{1,1}(D)}$. 

We now prove that if $|r|\, \theta$ belongs to $W^0_{1,1}(D)$, then $H\theta$ belongs to $W^0_{1}(D)$ and satisfies $\|H\theta\|_{W^0_1(D)}\leqs \||r|\, \theta\|_{W^0_{1,1}(D)}$. Theorem \ref{poincare} will then follow by interpolation.

Let $u$ be a non-vanishing vector field $u$. When we compute $H\theta  (z)[u(z)]$, we get
\begin{align}
H\theta  (z)[u(z)]&=\frac1{\vol(\Delta_n(\rho))}\int_{\over{\Lambda \in \Delta_n(\rho)}{t\in[0,1]}} \theta  (h_\Lambda(z,t))
 \left[\diffp{h_\Lambda}t (z,t),d_zh_\Lambda(z,t)[u]\right] \d t \d\Lambda.
\end{align}
Without restriction, we assume that the support of $\theta$ is included in a small neighborhood of $bD$. Therefore, in $H\theta$, we integrate only for $t\in[t_0,1]$, $t_0>0$. 
For $z\in D$ fixed, we decompose $[t_0,1]$ in 3 parts : $1-t\leq \frac\gamma2|r(z)|$, $1-t\geq \gamma |r(z)|$ and $\frac\gamma2|r(z)|\leq 1-t\leq \gamma |r(z)|$.
\subsection{Case $1-t\leq \frac\gamma2|r(z)|$}
We will use the following covering lemma :
\begin{lemma}\mlabel{covering}
 Let $K>0$ be arbitrary big and $\varepsilon_0>0$ be arbitrary small. If $c>0$ is small enough, there exists a sequence $(z_j)_{j\in\nn}$such that 
 \begin{enumerate}[(i)]
  \item \label{firstpoint} $D\setminus D_{-\varepsilon_0} \subset \bigcup_{j=0}^{+\infty} \cp{c|r(z_j)|}{z_j}$,
  \item \label{secondpoint} there exists $M$ such that all $z\in D\setminus D_{-\varepsilon_0}$, $z$ belongs to at most $M$ polydiscs  $\cp{cK|r(z_j)|}{z_j}$.
 \end{enumerate}
\end{lemma}
\pr The sequence is constructed as follows. Let $k$ be a non negative integer. We pick a point $z_{1}^{(k)}$ in  the boundary of $D_{-(1-c\kappa)^k\varepsilon_0}$ where $\kappa$ is a small positive number to be chosen later. We then pick up successively points $z_2^{(k)}, z_3^{(k)},\ldots$ in $bD_{-(1-c\kappa)^k\varepsilon_0}$ such that $\delta(z_j^{(k)},z_l^{(k)})\geq c\kappa (1-c\kappa)^k\varepsilon_0$ for all distinct $j$ and $l$. Then, there exists $\gamma>0$ such that for $j\neq l$, $\gamma\cp{c\kappa(1-c\kappa)^k\varepsilon_0}{z_j^{(k)}}\cap \gamma\cp{c\kappa(1-c\kappa)^k\varepsilon_0}{z_l^{(k)}}$ is empty and since $bD_{-(1-c\kappa)^k\varepsilon_0}$ is compact, this process stops at some rank $n_k$. Moreover, for all $z\in bD_{-(1-c\kappa)^k\varepsilon_0}$, there exists $j$ such that $z$ belongs to $\cp{c\kappa(1-c\kappa)^k\varepsilon_0}{z_j^{(k)}}$.

Let us prove that (\ref{firstpoint}) holds true. For $z\in D\setminus D_{-\varepsilon_0}$, let $k\in\nn$ be such that $(1-c\kappa)^{k+1}\varepsilon_0<|r(z)|\leq (1-c\kappa)^k\varepsilon_0$ and let $\lambda\in\rr$ be such that $\zeta=z+\lambda \eta_z $ belongs to $bD_{-(1-c\kappa)^k\varepsilon_0}$. 

Then there exists $j$ such that $\delta(\zeta, z_j^{(k)})\leq c\kappa(1-c\kappa)^k\varepsilon_0$. We also have $\delta(\zeta,z)=|\lambda|\leq c\kappa (1-c\kappa)^k\varepsilon_0$. Therefore $\delta (z,z^{(k)}_j) \leqs c\kappa \left|r\left(z_j^{(k)}\right)\right|$, and thus, if $\kappa$ has been chosen sufficiently small, $z$ belongs to $\cp{c\left|r\left(z_j^{(k)}\right)\right|}{z_j^{(k)}}$.

Now we prove (\ref{secondpoint}) of the lemma. 
Let $\zeta$ be a point in $D\setminus D_{-\varepsilon_0}$. If $\zeta$ belongs to $\cp{cK\left|r\left(z_j^{(k)}\right)\right|}{z_j}$, provided $c$ is small enough, we have $\frac12|r(\zeta)|\leq (1-c\kappa)^k\varepsilon_0\leq 2 |r(\zeta)|$. So there exist a finite number of $k$ such that $\zeta$ belongs to $\cp{cK\left|r\left(z_j^{(k)}\right)\right|}{z_j}$. For such a $k$, we put $$I_k=\left\{j\in\{1,\ldots, n_k\}\ / \ \zeta\in \cp{cK\left|r\left(z_j^{(k)}\right)\right|}{z_j}\right\}$$ and we show that $\#I_k$ is bounded, uniformly with respect to $k$.
We have for $C>0$, independant of $k$, $K$ and $c$, so big that $\cp {\frac cC(1-c\kappa)^k\varepsilon_0} {z_j^{(k)}}\cap \cp {\frac cC(1-c\kappa)^k\varepsilon_0} {z_l^{(k)}}=\emptyset$ for all $j\neq l$ :
\begin{align*}
 \sigma \left(\cup_{j\in I_k} \cp{cK\left|r\left(z_j^{(k)}\right)\right|}{z^{(k)}_j} \cap bD_{-(1-c\kappa)^k\varepsilon_0}\right)
 &\geq \sigma \left(\cup_{j\in I_k} \cp{\frac cCK\left|r\left(z_j^{(k)}\right)\right|}{z^{(k)}_j} \cap bD_{-(1-c\kappa)^k\varepsilon_0}\right)\\
 &\geq \sum_{j\in I_k}  \sigma \left(\cp{\frac cCK\left|r\left(z_j^{(k)}\right)\right|}{z^{(k)}_j} \cap bD_{-(1-c\kappa)^k\varepsilon_0}\right)\\
 &\geqs \sum_{j\in I_k}  \sigma \left(\cp{cK\left|r\left(z_j^{(k)}\right)\right|}{z^{(k)}_j} \cap bD_{-(1-c\kappa)^k\varepsilon_0}\right).
\end{align*}
Since $|r(\zeta)|\eqs \left|r\left(z_j^{(k)}\right)\right|$, we have $$\sigma \left(\cp{cK\left|r\left(z_j^{(k)}\right)\right|}{z^{(k)}_j} \cap bD_{-(1-c\kappa)^k\varepsilon_0}\right)\eqs \sigma \left(\cp{cK|r(\zeta)|}{\zeta} \cap bD_{-(1-c\kappa)^k\varepsilon_0}\right)$$ and so 
\begin{align*}
 \sigma \left(\cup_{j\in I_k} \cp{cK\left|r\left(z_j^{(k)}\right)\right|}{z^{(k)}_j} \cap bD_{-(1-c\kappa)^k\varepsilon_0}\right)
 &\geq \# I_k\cdot \sigma \left(\cp{cK|r(\zeta)|}{\zeta} \cap bD_{-(1-c\kappa)^k\varepsilon_0}\right).
 \end{align*}
On the other hand, since $\zeta$ belongs to $\cp{cK\left|r\left(z_j^{(k)}\right)\right|}{z^{(k)}_j}$ and since  $|r(\zeta)|\eqs \left|r\left(z_j^{(k)}\right)\right|$, there exists $C$ big such that 
 $\cp{cK\left|r\left(z_j^{(k)}\right)\right|}{z^{(k)}_j}\subset C\cp{cK|r(\zeta)|}{\zeta} $ and so
 \begin{align*}
 \sigma \left(\cup_{j\in I_k} \cp{cK\left|r\left(z_j^{(k)}\right)\right|}{z^{(k)}_j} \cap bD_{-(1-c\kappa)^k\varepsilon_0}\right)\leqs 
  \sigma \left(\cp{cK|r(\zeta)|}{\zeta} \cap bD_{-(1-c\kappa)^k\varepsilon_0}\right)
 \end{align*}
from which we get $\#I_k\leqs 1$.\qed

Now we proceed essentially as in \cite{Ale3}. Let $j$ be a non negative integer and set
\begin{align*}
\ij&:= \int_{\over{z\in\cp{c|r(z_j)|} {z_j}}{\over{\Lambda\in \Delta_n(\rho)}{t\in[1-\frac\gamma2\dz ,1]}} } \frac{|\theta(\hl)|\left[\diffp{h_\Lambda}t(z,t),\mathrm{d}_z\hl[u(z)]\right]}{k(z,u(z))\vol (\Delta_n(\rho))}
\mathrm{d}t\d\Lambda \d V(z).
\end{align*}
Lemma 2.17 from \cite{Ale3} implies that
\begin{align*}
\ij&\leqs  \int_{\over{z\in\cp{c|r(z_j)|} {z_j}}{\over{\Lambda\in \Delta_n(\rho)}{t\in[1-\frac\gamma2\dz ,1]}} } \frac{|\theta(\hl)|\left[\diffp{h_\Lambda}t(z,t),\mathrm{d}_z\hl[u(z)]\right]}{k\bigl(\hl,\d_z\hl[u(z)]\bigr)\cdot k\left(\hl,\diffp{h_\Lambda}t(z,t) \right)}
\mathrm{d}t\d\Lambda \d V(z).
\end{align*}
Then Proposition 2.12 from \cite{Ale3}  gives $\ij\leqs\sum_{k,l=1}^n \ijkl$ where
\begin{align*}
 \ijkl&:=
 \int_{\over{z\in\cp{c|r(z_j)|} {z_j}}{\over{\Lambda\in \Delta_n(\rho)}{t\in[1-\frac\gamma2\dz ,1]}} } \frac{|\theta(\hl)|[e_k(\hl),e_l(\hl)]}{k(\hl,e_k(\hl))\cdot k(\hl,e_l(\hl))}\d t\d\Lambda \d V(z).
\end{align*}
For fixed $z$ and $t$, we make the substitution $\zeta=\hl$, $\Lambda$ running over $\Delta_n(\rho)$. From \cite{Ale3} Lemma 2.15, when $|\Lambda|\leq \rho,$ the point  $\hl$ belongs to $C\frac{1-t}{\dz} \cp\dz z$ for some big $C>0$. Moreover, $\det_\rr \d_\Lambda h_\Lambda(z,t)\eqs \left(\frac{1-t}{\dz}\right)^{2n}(\det_\cc A(z))^2$ and Proposition 2.11 from \cite{Ale3} then gives 
$\det_\rr \d_\Lambda h_\Lambda(z,t)\eqs \left(\frac{1-t}{\dz}\right)^{2n} \vol(\cp\dz z)$. Therefore
\begin{align*}
 \ijkl&\leqs
 \int_{\over{z\in\cp{c|r(z_j)|} {z_j}}{\over{\zeta\in C\frac{1-t}{|r(z)|}\cp{|r(z)|}z}{t\in[1-\frac\gamma2\dz ,1]}} }
 \hskip-7pt \left(\frac{\dz}{1-t}\right)^{2n}\hskip-5pt  \frac{|\theta(\zeta)|[e_k(\zeta),e_l(\zeta)]}{\vol(\cp \dz z)\cdot \ k(\zeta,e_k(\zeta))\cdot k(\zeta,e_l(\zeta))}\d V(\zeta)\d t \d V(z).
\end{align*}
Now we want to change the order of integration. Intuitively, since $t$ is sufficiently close to $1$, $\zeta$ will belong to $\cp{cK|r(z_j)|}{z_j}$ for some $K$ and $z$ will belong to  
$\cp{|r(\zeta)|}\zeta$.

If $\gamma$ is small enough, $|r(z)|\eqs |r(\zeta)|$ because $\frac{1-t}{|r(z)|}\leq \frac\gamma2$ and because $\zeta$ belongs to $C\frac{1-t}{|r(z)|}\cp{|r(z)|}z $. This implies, if $\gamma$ is even smaller, that $t$ belongs to $[1-\frac12|r(\zeta)|,1]$. We also have $\delta(\zeta,z_j)\leqs \delta(\zeta,z)+\delta(z,z_j)\leqs c|r(z_j)|$ so $\zeta$ belongs to $\cp{cK|r(z_j)|}{z_j}$ for some $K$ which does not depend on $j$, $c$ or $\zeta$. Therefore, $c$ can be chosen small so that Lemma \ref{covering} holds true.

Since $\zeta$ belongs to $C\frac{1-t}{\dz}\cp \dz z$, we can write $\zeta=z+C\frac{1-t}{\dz}\mu v$ with $\mu\in\cc$, $v\in\cc^n$, $|v|=1$, such that $|\mu|<\tau(z,v,\dz)$. Provided $\gamma$ is small enough, we have $\dz\eqs\dzeta$ and $\tau(z,v,\dz)\eqs\tau(\zeta,v,\dzeta)$. Therefore $z=\zeta-\frac{1-t}{\dz}\mu v$ with $|\mu|\leqs \tau(\zeta,v,\dzeta)$ and there exists $K'>0$ big, such that $z$ belongs to $K'\frac{1-t}{\dzeta}\cp{\dzeta}{\zeta}$.\\
Therefore, the set 
$\left\{(z,t,\zeta),\ z\in\cp {|r(z_j)|}{z_j} ,\ t\in[1-\frac\gamma2\dz ,1],\ \zeta\in C\frac{1-t}{\dz } \cp\dz z\right\}$ is included in 
$\left\{(z,t,\zeta),\ \zeta\in \cp {cK|r(z_j)|}{z_j},\ t\in[1-\frac12\dzeta,1],\ z\in K'\frac{1-t}{\dzeta} \cp{\dzeta}\zeta\right\}$.\\
Moreover, $\vol(\cp\dz z)\eqs\vol( \cp {\dzeta}\zeta)$ which gives
\begin{align*}
 \ijkl&\leqs
 \int_{\over{\zeta\in \cp {cK|r(z_j)|}{z_j}}{\over{t\in[1-\frac12\dzeta,1]}{z\in K'\frac{1-t}{\dzeta}\cp{\dzeta}\zeta}} }
 \hskip-9pt\left(\frac{|r(\zeta)|}{1-t}\right)^{2n}\hskip-7pt \frac{|\theta(\zeta)|[e_k(\zeta),e_l(\zeta)]}{\vol(\cp \dzeta \zeta)\cdot \ k(\zeta,e_k(\zeta))\cdot k(\zeta,e_l(\zeta))}\d V(\zeta)\d t \d V(z).
\end{align*}
We integrate successively with respect to $z$ and $t$ and get 
\begin{align*}
 \ijkl&\leqs
 \int_{{\zeta\in \cp {cK|r(z_j)|}{z_j}}}
  \frac{|r(\zeta)||\theta(\zeta)|[e_k(\zeta),e_l(\zeta)]}{\ k(\zeta,e_k(\zeta))\cdot k(\zeta,e_l(\zeta))}\d V(\zeta).
\end{align*}
Now, summing over $j\in\nn$, we get, since any $\zeta$ belongs to at most $M$ polydiscs $\cp {cK|r(z_j)|}{z_j}$
\begin{align*}
 &\int_{\over{z\in D}{\over{\Lambda\in \Delta_n(\rho)}{t\in[1-\frac\gamma2\dz ,1]}}} \frac{|\theta(\hl)|\left[\diffp{h_\Lambda}t(z,t),\mathrm{d}_z\hl[u(z)]\right]}{k(z,u(z))\vol (\Delta_n(\rho))}
\mathrm{d}t\d\Lambda \d V(z)
\\
&\leqs \sum_{j=0}^{+\infty} \sum_{k,l=1}^n \int_{{\zeta\in \cp {cK|r(z_j)|}{z_j}}}
  \frac{|r(\zeta)||\theta(\zeta)|[e_k(\zeta),e_l(\zeta)]}{\ k(\zeta,e_k(\zeta))\cdot k(\zeta,e_l(\zeta))}\d V(\zeta)\\
&\leqs  \sum_{k,l=1}^n \int_{\zeta\in D}
  \frac{|r(\zeta)||\theta(\zeta)|[e_k(\zeta),e_l(\zeta)]}{ \ k(\zeta,e_k(\zeta))\cdot k(\zeta,e_l(\zeta))}\d V(\zeta)\\
&\leqs \||r| \theta\|_{W^0_{1,1}(D)}.
\end{align*}

\subsection{Case $1-t\geq \gamma |r(z)|$}
Here we want to estimate 
\begin{align*}
 \ii&:=\int_{\over{z\in D}{\over{\Lambda\in \Delta_n(\rho)}{t\in[t_0,1-\gamma\dz ]}}} \frac{|\theta(\hl)|\left[\diffp{h_\Lambda}t(z,t),\mathrm{d}_z\hl[u(z)]\right]}{k(z,u(z))\vol (\Delta_n(\rho))}
\mathrm{d}t\d\Lambda \d V(z).
\end{align*}
Lemma 2.20 of \cite{Ale3} implies
\begin{align*}
 \ii&\leqs\int_{\over{z\in D}{\over{\Lambda\in \Delta_n(\rho)}{t\in[t_0,1-\gamma\dz ]}}} 
 \left(\frac{1-t}{|r(z)|}\right)^{1-\frac1m} \frac{|\theta(\hl)|\left[\diffp{h_\Lambda}t(z,t),\mathrm{d}_z\hl[u(z)]\right]\d t\d\Lambda \d V(z)}
 {k\bigl(\hl,\d_z\hl[u(z)]\bigr)\cdot k\left(\hl,\diffp{h_\Lambda}t(z,t) \right)}.
\end{align*}
and Proposition 2.12 of \cite{Ale3} gives $\ii\leqs \sum_{j,k=1}^n \iijk$ where
\begin{align*}
 \iijk&:=
 \int_{\over{z\in D}{\over{\Lambda\in \Delta_n(\rho)}{t\in[t_0,1-\gamma\dz ]}}} 
 \left(\frac{1-t}{|r(z)|}\right)^{1-\frac1m} \frac{|\theta(\hl)|\left[e_j(\hl),e_k(\hl)\right]\d t\d\Lambda \d V(z)}
 {k\bigl(\hl,e_j(\hl)\bigr)\cdot k\left(\hl,e_k(\hl)\right)}.
\end{align*}
We now make the substitution $\zeta=\hl$. Since, Proposition 2.11 of \cite{Ale3}, $\det_\rr \hl=\det_\rr A(tz)\eqs \vol\left(\cp {|r(tz)|}{tz}\right)$ and since, Lemma 2.18 of \cite{Ale3}, $\{\hl\ / \ \Lambda \in \Delta_n(\rho)\}$ is a subset of $Ct\rho \cp {|r(tz)|}{tz}$, we have
\begin{align*}
 \iijk\hskip -2pt &\leqs \hskip -2pt
 \int_{\over{z\in D}{\over{t\in[t_0,1-\gamma\dz ]}{\zeta\in Ct\rho \cp{|r(tz)|}{tz}}}} \hskip -8pt 
 \left(\frac{1-t}{|r(z)|}\right)^{1-\frac1m}\hskip -8pt \frac{|\theta(\zeta)|\left[e_j(\zeta),e_k(\zeta)\right]}
 {k\bigl(\zeta,e_j(\zeta)\bigr)\cdot k\left(\zeta,e_k(\zeta)\right) \vol\left({\cp{|r(tz)|}{tz}}\right)}\d t\d V(\zeta) \d V(z).
\end{align*}
We want to apply Fubini's theorem. For $\zeta\in D$ fixed, if t and $z$ are such that $\zeta$ belongs to $Ct\rho\cp{|r(tz)|}{tz}$, then $|r(tz)|\eqs |r(\zeta)|$ if $\rho $ is small enough and since, Corollary 2.19 of \cite{Ale3},   $|r(tz)|\eqs 1-t$, we have $1-t\eqs |r(\zeta)|$. Moreover
\begin{align*}
 \delta(\zeta,z)&\leqs \delta(\zeta,tz)+\delta(tz,z)\\
 &\leqs |r(tz)|+|tz+z|\\
 &\leqs |r(\zeta)|+1-t \leqs |r(\zeta)|.
 \end{align*}
So $z$ belongs to $\cp {K|r(\zeta)|}\zeta$ for some big $K$. Finally, if $\rho$ is small enough,  since $\zeta$ belongs to $Ct\rho\cp {|r(tz)|}{tz}$, $\vol \cp{|r(tz)|} {tz}\eqs \vol \cp {|r(\zeta)|}\zeta$. We thus have
\begin{align*}
  \iijk&\leqs
 \int_{\over{\zeta\in D}{\over{|r(\zeta)|\leqs 1-t \leqs |r(\zeta)|}{z \in D\cap \cp{K|r(\zeta)|}\zeta }}} 
 \left(\frac{1-t}{|r(z)|}\right)^{1-\frac1m} \frac{|\theta(\zeta)|\left[e_j(\zeta),e_k(\zeta)\right]}
 {k\bigl(\zeta,e_j(\zeta)\bigr)\cdot k\left(\zeta,e_k(\zeta)\right) \vol{\cp{|r(\zeta)|}{\zeta}}}\d t\d V(\zeta) \d V(z).
\end{align*}
Now, using  the 2 inequalities $\int_{z\in D\cap \cp {K|r(\zeta)|}\zeta} |r(z)|^{\frac1m-1}d\lambda(z)\leqs |r(\zeta)|^{\frac1m-1} \vol \cp {|r(\zeta)|}\zeta$ and $\int_{|r(\zeta)|\leqs 1-t\leqs |r(\zeta)|} (1-t)^{1-\frac1m}\d t \leqs |r(\zeta)|^{2-\frac1m}$, we get
\begin{align*}
  \iijk&\leqs
 \int_{\zeta\in D} 
\frac{|r(\zeta)| |\theta(\zeta)|\left[e_j(\zeta),e_k(\zeta)\right]}
 {k\bigl(\zeta,e_j(\zeta)\bigr)\cdot k\left(\zeta,e_k(\zeta)\right)} \d V(\zeta),
\end{align*}
from which we conclude that $\ii\leqs \||r|\,\theta\|_{W^0_{1,1}(D)}$.
\subsection{Case $\frac\gamma2|r(z)|\leq 1-t\leq \gamma |r(z)|$}
For $j\in\nn$, we set
\begin{align*}
\iiij:= \int_{\over{z\in\cp{c|r(z_j)|} {z_j}\cap D}{\over{t\in[1-\gamma\dz ,1-\frac\gamma2\dz]} {\Lambda\in \Delta_n(\rho)}}} \frac{|\theta(\hl)|\left[\diffp{h_\Lambda}t(z,t),d_z\hl[u(z)]\right]}{k(z,u(z))\vol (\Delta_n(\rho))}
\d\Lambda \d t\d V(z).
\end{align*}
Combining Lemma 2.21 and Lemma 2.12 of \cite{Ale3} gives  $\iiij\leqs\sum_{k,l=1}^n \iiijkl$ where 
\begin{align*}
 \iiijkl&:=\hskip -2pt
 \int_{\over{z\in\cp{c|r(z_j)|} {z_j}}{\over{\Lambda\in \Delta_n(\rho)}{t\in[1-\gamma\dz, 1-\frac\gamma2\dz }} }\hskip -5pt \frac{|\theta(\hl)|[e_l(\hl),e_k(\hl)]}{k(\hl,e_l(\hl))\cdot k(\hl,e_k(\hl))}\d t\d\Lambda \d V(z).
\end{align*}

Now we make the substitution $\zeta=\hl$, $\Lambda\in \Delta_n(\rho)$. By Lemma 2.21 of \cite{Ale3}, $\hl$ belongs to $\cp{c|r(z)|}z$ and $|r(\hl)|\eqs \dz$ if $\gamma$ is small enough. As in \cite{Ale3}, Subsection 2.5, $\det_\rr (\d_\Lambda \hl)\geqs \vol (\cp{\dz}z)$, thus
\begin{align*}
 \iiijkl
 &\leqs \int_{\over{z\in\cp{c|r(z_j)|} {z_j}}{\over{\zeta\in\cp{c\dz}z}{t\in[1-\gamma\dz, 1-\frac\gamma2\dz }} } \frac{|\theta(\zeta)|[e_l(\zeta),e_k(\zeta)]}{k(\zeta,e_l(\zeta))\cdot k(\zeta,e_k(\zeta))\vol \cp{\dz}z}\d t\d V(\zeta) \d V(z)\\
  &\leqs \int_{\over{z\in\cp{c|r(z_j)|} {z_j}}{\zeta\in\cp{c\dz}z}} \frac{|r(z)||\theta(\zeta)|[e_l(\zeta),e_k(\zeta)]}{k(\zeta,e_l(\zeta))\cdot k(\zeta,e_k(\zeta))\vol \cp{\dz}z}\d V(\zeta) \d V(z)
\end{align*}
For $\zeta$ in $\cp {c\dz}z$ and $z$ in $\cp{c|r(z_j)|}{z_j}$, if $c$ is small enough, $|r(\zeta)|\eqs \dz\eqs |r(z_j)|$  $\vol\cp{\dz}z\eqs\vol\cp{|r(\zeta)|}\zeta$ and $\delta(\zeta,z_j)\leqs\delta(\zeta,z)+\delta(z,z_j)\leqs c|r(z_j)|$, so $\zeta$ belongs to $\cp{cK|r(z_j)|}{z_j}$ for some big $K$, not depending on $c$, $\zeta,$ $z$ or $z_j$. The point $z$ also belongs to $\cp{cK|r(\zeta)|}\zeta$ if $K$ is big enough so
\begin{align*}
 \iiijkl
 &\leqs \int_{\over{\zeta\in\cp{cK|r(z_j)|} {z_j}}{z\in\cp{cK|r(\zeta)|}\zeta}} \frac{|r(\zeta)|\,|\theta(\zeta)|[e_l(\zeta),e_k(\zeta)]}{k(\zeta,e_l(\zeta))\cdot k(\zeta,e_k(\zeta))\vol \cp{|r(\zeta)|}\zeta}\d V(\zeta) \d V(z)\\
 &\leqs \int_{{\zeta\in\cp{cK|r(z_j)|} {z_j}}} \frac{|r(\zeta)|\,|\theta(\zeta)|[e_l(\zeta),e_k(\zeta)]}{k(\zeta,e_l(\zeta))\cdot k(\zeta,e_k(\zeta))}\d V(\zeta).\\
\end{align*}
Since, Lemma \ref{covering}, any $\zeta$ belongs to at most $M$ polydiscs $\cp {cK|r(z_j)|}{z_j}$, we get
\begin{align*}
\int_{\over{z\in D}{\over{t\in[1-\gamma\dz ,1-\frac\gamma2\dz]} {\Lambda\in \Delta_n(\rho)}}} \frac{|\theta(\hl)|\left[\diffp{h_\Lambda}t(z,t),d_z\hl[u(z)]\right]}{k(z,u(z))\vol (\Delta_n(\rho))}
\d\Lambda \d t\d V(z)&\leqs\||r|\,\theta\|_{W^0_{1,1}(D)},
\end{align*}
 which conclude the proof of Theorem \ref{poincare}.
\bibliographystyle{siam}
\bibliography{../bibliographie/bibliographie}
\end{document}